\documentclass[10pt]{article}
\usepackage{amssymb}
\usepackage{epsfig}
\usepackage{graphics}
\input amssym.def
\input amssym
\newtheorem{theorem}{Theorem}[section]
\newtheorem{corollary}{Corollary}[section]
\newtheorem{lemma}{Lemma}[section]

\newtheorem{remark}{Remark}[section]
\newtheorem{definition}{Definition}[section]
\newtheorem{proposition}{Proposition}[section]

\newcommand{\nb}{\mbox{\tiny ${\bullet\atop\bullet}$}}

\newcommand{\hD}{\hat{\mathcal{D}}}

\newcommand{\bea}{\begin{eqnarray}}
\newcommand{\halmos}{\rule{1ex}{1.4ex}}
\newcommand{\ds}{\displaystyle}
\newcommand{\epfv}{\hspace*{\fill}\mbox{$\halmos$}\vspace{1em}}

\newcommand{\eea}{\end{eqnarray}}
\newcommand{\nn}{\nonumber \\}
\newcommand{\be}{\begin {equation}}

\newcommand{\ee}{\end{equation}}

\title{{Formal differential operators, vertex
operator algebras and zeta--values , II}}
\author{Antun Milas \\
Department of Mathematics \\ University of Arizona,
Tucson, AZ 85721}
\date{}
\begin{document}
\maketitle
\small{
\begin{abstract}
We introduce certain correlation functions (graded $q$--traces)
associated to vertex operator algebras and superalgebras which we
refer to as $n$--point functions. These naturally arise in the studies
of representations of Lie algebras of 
differential operators on the circle \cite{Le1}--\cite{Le2}, \cite{M}.
We investigate their properties and consider the corresponding
graded $q$--traces in parallel with the passage from genus $0$ to
genus $1$ conformal field theory. By using the vertex operator
algebra theory we analyze in detail correlation 
functions in some particular cases. We obtain elliptic
transformation properties for $q$--traces and the corresponding
$q$--difference equations. In particular, our construction 
leads to certain correlation functions and $q$--difference equations investigated by S. Bloch and A. Okounkov \cite{BO}.
\end{abstract}

\renewcommand{\theequation}{\thesection.\arabic{equation}}
\setcounter{equation}{0}

\section{Introduction}
This is a continuation of \cite{M}. In this part we study certain correlation functions built up from
the iterates of vertex operators introduced in Part I \cite{M}.

Let $V$ be an arbitrary vertex operator (super)algebra and
$M$ a $V$--module. Sometimes we will weakened this property by assuming that
$M$ is only a subspace of a $V$--module invariant with respect to
certain operators.
Let $u_i \in V$, $i=1,...,n$, $w \in M$ and $w' \in M'$, where
$M'$ is the (restricted) dual space of $M$ as defined in \cite{FHL}.
Suppose that $x_i$'s are commuting formal variables as in Part I \cite{M}.
In the vertex operator algebra theory one usually studies the following
formal {\em matrix coefficients}
\begin{equation} \label{000}
\langle w', Y(u_1,x_1) \cdots Y(u_n,x_n)w \rangle ,
\end{equation}
and the corresponding graded $q$--traces (cf. \cite{Zh1})
\begin{equation} \label{00a}
{\rm tr}|_{M} X(u_1,x_1) \cdots X(u_n,x_n) q^{\bar{L}(0)}.
\end{equation}
Matrix coefficients of the type (\ref{000})
are studied in \cite{FLM} in both the formal and
the analytic context (meaning that $x_i$ are set to be
complex variables). In particular in  \cite{Hu} expressions of the form
(\ref{000}) were used for construction of the genus zero
meromorphic conformal field theory.
On the other hand the graded traces function of the form (\ref{00a})
are related to the genus one meromorphic conformal
field theory. More precisely, the expression of the form
(\ref{00a}) give a vector in a genus--one conformal block associated
to $V$ \cite{Zh1}.

In Part I \cite{M} we have studied a relationship between various
normal orderings and classical Lie subalgebras of the Lie algebra of
(super)differential operators on the circle. The most interesting
representations of $\hD$ (unitary representations for instance)
can be constructed by using certain quadratic operators in terms
of free fields \cite{AFOQ},\cite{FKRW}, \cite{KWY}.
As we already noticed, these quadratic
operators are closely related to iterates $X(Y[u,y]v,x)$ (cf.
\cite{Le1}--\cite{Le3}, \cite{M}). Therefore it is very natural to consider matrix
coefficients of the form
\begin{equation} \label{02}
\langle w', X(Y[u_1,y_1]v_1,x_1) \cdots
X(Y[u_n,y_n]v_n,x_n)w
\rangle,
\end{equation}
which we will call iterated $2n$--point function
\footnote{$2n$ refers to the number of formal (or
complex) variables}.
Also we consider a
$2n$--point function :
\begin{equation} \label{01}
\langle u'_{n+1}, X(u_1,x_1 t_1)X(v_1,x_1) \cdots X(u_n,x_nt_n)
X(v_n,x_n) u_{n+1} \rangle.
\end{equation}
The first result in our paper is Theorem \ref{main1} which relates
the matrix
coefficients (\ref{02}) and (\ref{01}).
In parallel with (\ref{000}) and (\ref{00a}) it is natural to consider
the corresponding graded $q$--traces:
\bea \label{002}
{\rm tr}|_M X(Y[u_1,y_1]v_1,x_1) \cdots
X(Y[u_n,y_n]v_n,x_n)q^{\bar{L}(0)},
\eea
Here $\bar{L}(0)=L(0)-\frac{c}{24}$.
Again, (cf. Proposition \ref{main2}) (\ref{002}) is closely related to
\bea \label{001}
{\rm tr}|_M X(u_1,x_1 t_1)X(v_1,x_1) \cdots X(u_n,x_nt_n)
X(v_n,x_n)q^{\bar{L}(0)}.
\eea


Besides the $q$--traces that we already mentioned, one would like to consider
(for reasons that will become clear in a moment)
an expansion of the $q$--trace (\ref{001}) in powers of the $x_i$'s:
$$ {\rm tr}|_M X(u_1,x_1 t_1)X(v_1,x_1) \cdots X(u_n,x_nt_n)
Y(v_n,x_n)q^{\bar{L}(0)}=$$
$$\sum_{\alpha} ({{\rm tr}|_{M}X(u_1,x_1 t_1)X(v_1,x_1) \cdots X(u_n,x_nt_n)
X(v_n,x_n)q^{L(0)}})_{\alpha} x^{\alpha}$$
(here we use the multi--index notation $x^{\alpha}=x_1^{\alpha_1} \cdots
x_n^{\alpha_n}$ ).
Especially interesting is the {\em constant term} with respect to
the $x$--variables, i.e.,
\bea \label{003}
&&{\rm Coeff}_{x_1^0 \ldots x_n^0}
{\rm tr}|_{M} X(u_1,x_1 t_1)X(v_1,x_1) \cdots X(u_n,x_nt_n)
X(v_n,x_n)q^{\bar{L}(0)}=\nn
&& {\rm tr}|_{M} o(X(u_1,x_1 t_1)X(v_1,x_1)) \cdots o(X(u_n,x_nt_n)
X(v_n,x_n))q^{\bar{L}(0)}
\eea
where $o(a)=a({\rm wt}(a)-1)$ (cf. \cite{Zh1}).
This new formal expression
depends on $n$--parameters: $t_1$,...,$t_n$.
Also we are interested in
\bea \label{004}
&&{\rm Coeff}_{x_1^0 \ldots x_n^0}
{\rm tr}|_{M} X(Y[u_1,y_1]v_1,x_1) \cdots X(Y[u_n,y_n]v_n,x_n)
q^{\bar{L}(0)}=\nn
&& {\rm tr}|_{M} o(Y[u_1,y_1]v_1) \cdots o(Y[u_n,y_n]v_n)q^{\bar{L}(0)}.
\eea
Any expression of the form (\ref{003}) or (\ref{004}) (possibly normalized by
${\rm tr}|_M q^{L(0)}$) we shall refer to as {\em $n$--point correlation
function} (or simply, $n$--point function).
This term is widely used in statistical physics and random matrix
theory where an $n$--point correlation function depending on
$y_1,...,y_n$ (or $t_1$,...,$t_n$)
can be defined, for example, as an
integral (usually normalized by the partition function)
of a probability density function depending on integration
variables $x_1,...,x_n$.
Because extracting the zero--th term in
(\ref{003}) resembles (complex) integration in variables $x_1$,...,$x_n$
we decided to use this terminology (cf. \cite{BO}).


In our applications $M$ will be a $V$--module stable
under the Fourier modes of $X(u_i,x t) X(v_i,x)$.
In the case of rational vertex operator algebras with some 
additional properties (cf. \cite{Zh1}), the $q$--trace (\ref{00a})
rise to a doubly periodic functions (in this case even more is true; the vector space of
all characters of $V$--modules is invariant under $SL(2, \mathbb{Z})$ action).
However, to prove ellipticity one needs much weaker assumptions.

Here is a short overview:
\begin{itemize}
\item[(i)]
In Section 2 we derive a precise relationship
between (\ref{02}) and (\ref{01}). We also
relate the corresponding $q$--traces (\ref{002}) and (\ref{001})
and the corresponding $n$--point functions (\ref{003}) and (\ref{004}).

\item[(ii)] In Section 3 we apply our results from Section 2 in
the case of the infinite--wedge vertex operator (super)algebra
$\mathcal{F}$ and its charge $m$ subspace $\mathcal{F}_m$. 
It is well--known that these subspaces are $\hat{\mathcal{D}}$--modules. 
We compute $q$--traces (\ref{002}) and (\ref{001}) in the most
interesting case (when the vectors $u_i$ and $v_i$ are chosen to
be the fermionic generators) and the corresponding $n$--point functions.
In particular, the $n$--point correlation functions (\ref{004}) rise to the
$n$--point functions introduced by Bloch and Okounkov \cite{BO}, \cite{O}.
The rest of Section 3 is devoted to studies of the $q$--difference
equations satisfied by these $n$--point functions. We obtain
several explicit formulas for $1$-- and $2$--point functions by
using different methods.

\item[(iii)]
In Section 4, along the lines of Section 3, we consider a fermionic
vertex operator superalgebra and the corresponding $n$--point functions
associated with it. This case rises to a representation of an {\em orthogonal}
Lie algebra of differential operators on the circle $\hat{\mathcal{D}}^-$.

\item[(iv)] In Section 5, we consider a vertex operator algebra
associated to a free boson and the corresponding $n$--point functions
associated to the bosonic generators. This case was studied in
\cite{M} in connection with a {\em symplectic} Lie algebra of differential
operators on the circle $\hat{\mathcal{D}}^+$ \cite{Bl}, \cite{M} and
the zeta--regularization procedure.

\item[(v)] In Appendix A we prove certain elementary lemmas
necessary to deal with $q$--difference equations in sections 3,4 and
5. Finally, in Appendix B we gave a different proof of the
so--called ``recursion formula'' for the $q$--graded traces
originally due to Zhu \cite{Zh1}.
\end{itemize}
{\em n.b.} This paper was originally written in 2000 and it is a part
of authors Ph.D. thesis. In the meantime we have noticed that 
Okounkov \cite{O}, Miyamoto \cite{Mi}
and Dong, Mason and Nagatomo \cite{DMN} 
obtained some interesting results related to our paper.

\renewcommand{\theequation}{\thesection.\arabic{equation}}
\setcounter{equation}{0}

\section{Correlation functions}

As in Part I \cite{M} we will use formal calculus as developed in
\cite{FLM}. We denote by $x$, $y$, $t$, $x_i$, $t_i$ and $y_i$, etc.,
commuting formal variables and we take the liberty of using the same notation when formal variables are
replaced by complex numbers. From the context it should be clear
whether the variables are formal or complex.

\subsection{Normal ordering procedure} \label{trivial}
Let $(V,Y,{\bf 1},\omega)$ be a vertex operator algebra ( see
\cite{FLM} for the definition ) and $u,v \in V$. From the Jacobi
identity for vertex operator algebras it follows that \be
\label{lll} Y(u(-1)v,x)=\nb Y(u,x)Y(v,x) \nb:=
Y^-(u,x)Y(v,x)+Y(v,x)Y^+(u,x), \ee where $Y^+(v,x)=\ds{\sum_{n
\geq 0} v(n) x^{-n-1}}$ is the singular and
$Y^-(v,x)=\ds{\sum_{n<0} v(n) x^{-n-1}}$ is the regular part of
$Y(v,x)$. Also, we define
$$\nb Y(u,x)Y(v,y) \nb= Y^-(u,x)Y(v,y)+Y(v,y)Y^+(u,x).$$
Let us recall the definition of the $X$--operator as in the
first part \cite{M}:
$$X(u,x)=Y(x^{L(0)}u,x),$$
where $u$ is a homogeneous vector.
We extend this definition by the linearity for every $u \in V$.
The $X$--operator also admits a splitting into the regular and singular part, but then
$X^-(u,x) \neq Y^-(x^{L(0)}u,x)$.
Thus we have two different normal orderings for $X$ operators
$$\nb' X(u,x)X(v,x) \nb'=X^-(u,x)X(v,x)+X(v,x)X^+(u,x)$$
and
$$\nb  X(u,x)X(v,x) \nb =
Y^-(x^{L(0)}u,x)X(v,x)+X(v,x)Y^+(x^{L(0)}u,x).$$ It is more
convenient to work with $\nb \ \nb$ instead of $\nb' \ \nb'$. If
we suppose that $$u, v \in V_1 \ \ {\rm and} \ \ [v(0),u(n)]=0,$$
for every $n$, then \bea \label{nobf} &&  \nb' X(u,x)X(v,x)
\nb'=\left\{\sum_{n \leq 0} u(n) x^{-n} \right\}X(v,x)+ X(v,x)
\left\{\sum_{n > 0} u(n) x^{-n} \right\}= \nn &&= \left\{\sum_{n <
0} u(n) x^{-n} \right\} X(v,x)+u(0)X(v,x)+ X(v,x) \left\{\sum_{n
> 0} u(n) x^{-n} \right\} \nn
&&= \left\{ x \sum_{n < 0} u(n) x^{-n-1} \right\}X(v,x)+X(v,x) x
\left\{\sum_{n \geq  0} u(n) x^{-n-1} \right\} \nn
&&=xY^-(u,x)X(v,x)+X(v,x)xY^+(u,x)=\nb X(u,x)X(v,x) \nb, \eea
meaning that there is no ambiguity which normal ordering we are
using. This case arises when we study bosons and fermions
\footnote{ All results about normal ordering hold with minor
modifications for the vertex operator superalgebras. In
particular, the formula (\ref{lll}) has to be replaced by $\nb
Y(u,x)Y(v,x) \nb:=
Y^-(u,x)Y(v,x)+(-1)^{p(u)p(v)}Y(v,x)Y^+(u,x)$.}.
%
\begin{proposition}
Suppose that \be \label{condition} u(n)v=c_{u,v}\delta_{n,{\rm
wt}(u)+{\rm wt}(v)-1}{\bf 1}, \ee where $c_{u,v} \in \mathbb{C}$.
Then
\begin{itemize}
\item[(a)]
$$ \nb X(u,x_1)X(v,x_2) \nb=\nb X(v,x_2)X(u,x_1) \nb $$
and
\item[(b)]
$$ \nb ' X(u,x_1)X(v,x_2) \nb '=\nb X(u,x_1)X(v,x_2) \nb.$$
\end{itemize}
\end{proposition}
We will not use this result in the rest of the paper so we leave the
proof to the reader.

\subsection{Iterated $2n$--point  functions }

In what follows we shall always use the following binomial
expansion conventions. An expression $\frac{1}{(x-y)^k}$, where
$x$ and $y$ are formal variables, is understood to be expanded (by
the binomial theorem) in non--negative powers of $y$. Note that
the order of variables is important. Also we allow $x$ or $y$ (but
not both) to be a complex number.

As in \cite{FHL} we will be considering formal series
of the form
\begin{equation} \label{genformcor}
\frac{g(x_1,...,x_n)}{\displaystyle{\prod_{i=1}^n x_i^{r_i}
\prod_{j < k} (x_j -x_k)^{s_{jk}}}} \in \mathbb{C}[[x_1^{\pm
1},...,x_n^{\pm 1}]],
\end{equation}
where $g(x_1,...,x_n) \in {\mathbb C}[x_1,...,x_n]$.

It is well--known (cf. \cite{FHL}, \cite{FLM}) that
$$\langle u'_{n+1}, Y(u_1,x_1) \cdots Y(u_n,x_n)u_{n+1} \rangle=
f(x_1,...,x_n),$$ for some $f(x_1,...,x_n)$ of the form
(\ref{genformcor}). After we replace the formal variables with
complex variables (\ref{genformcor}) converges to a rational
function inside the domain $|x_1|> \cdots >|x_n|>0$.

Let $p_i, s_j, p_{i,j}, r_{i,j},s_{i,j}$ and $u_{i,j}$ be natural numbers and
$g(x_i,t_j)_{i,j=1}^n \in {\mathbb C}[x_i,t_j]_{i,j=1}^n$.
Then
\bea \label{genform}
&& f(x_i,t_j)_{i,j=1}^n= \\
&& = \frac{g(x_i,t_j)_{i,j=1}^n}{\displaystyle{ \prod_{i=1}^n
(t_i-1)^{p_i} t_i^{r_i} \prod_{j=1}^n x_j^{s_j} \prod_{i
<j}}(t_ix_i-t_jx_j)^{p_{i,j}}(t_ix_i-x_j)^{r_{i,j}}(x_i-t_jx_j)^{s_{i,j}}
(x_i-x_j)^{u_{i,j}}} \nonumber \eea is a well--defined element of
$\mathbb{C}[[x_1^{\pm 1},...,x_n^{\pm 1},t^{\pm 1}_1,...,t^{\pm
1}_n]]$ where  we expand the denominator by using binomial
expansion. Viewed as a rational function (\ref{genform}) has poles
on the divisors: $x_i=x_j$, $x_i=0$, $t_i=0$, $t_i=1$,
$x_i=\frac{t_j x_j}{t_i}$ , $x_i=t_jx_j$ and $x_i=x_jt_i^{-1}$. It
is also clear that \bea \label{genform1} \langle u'_{n+1},
X(u_1,x_1 t_1)X(v_1,x_1) \cdots X(u_n,x_nt_n) X(v_n,x_n) u_{n+1}
\rangle= g(x_i,t_j)_{i,=1}^n, \eea where
$f(x_1,...,x_n,t_1,...,t_n)$ is of the form (\ref{genform}). We
may assume $t_{i_1}=...=t_{i_n}$ for some $1 \leq i_1 \leq \ldots
\leq i_k \leq n$; in this case (\ref{genform}) is still
well--defined.

It is important to notice that in (\ref{genform}) we cannot perform
the substitution $t_i=e^{y_i}$, $i=1,...,n$,
since the binomial expansion convention  is not applicable to the
expansion  $\frac{1}{(e^y-1)^k}$.
Rather we use a different convention.
As before $\frac{1}{(e^{y_i}x_i-e^{y_j}x_j)^k}$ are expanded
by using binomial expansion, and
factors of the form
$$\frac{1}{(e^{x}-1)^k},$$
$k \in \mathbb{N}$
stand for the formal multiplicative inverse of $(e^{x}-1)^k$, i.e.
$$\frac{x^{-k}}{(1+(\frac{x}{2!}+\frac{x^2}{3!}+...))^k} \in \mathbb{C}((x)),$$
where we use binomial expansion with respect to
$\frac{x}{2!}+\frac{x^2}{3!}+ \cdots$. Sometimes we will view
$\frac{1}{(e^{y}-1)^k}$ as a meromorphic function in $y$.
With these conventions we show
\begin{proposition} \label{welld}
$f(x_i,e^{y_j})_{i=1,j=1}^n$ is a well--defined element of
$$\mathbb{C}((y_1,...,y_n))[[x_1^{\pm 1},...,x_n^{\pm 1}]].$$
\end{proposition}
{\em Proof:} It is enough to show that in the expansion
$$ f(x_i,e^{y_j})_{i=1,j=1}^n=\sum_{\alpha} x^{\alpha} f_{\alpha}(e^{y_i})_{i=1}^n,$$
(we use a multi index notation $\alpha=(\alpha_1,...,\alpha_n)$) for each $\alpha$,
$f_{\alpha}(e^{y_i})_{i=1}^n$ is a well--defined element of
$\mathbb{C}((y_1,...,y_n))$.
\bea \label{induction}
&& {\rm Coeff}_{x^\alpha}
 \frac{g(x_i,e^{y_j})_{i,j=1}^n}
{\displaystyle{\prod_{i<j}}(e^{y_i}x_i-e^{y_j}x_j)^{p_{i,j}}
(e^{y_i}x_i-x_j)^{r_{i,j}}(x_i-e^{y_j}x_j)^{s_{i,j}}
(x_i-x_j)^{u_{i,j}}}=\nn
&&= {\rm Coeff}_{x_1^{\alpha_1}} ( {\rm Coeff}_{x_2^{\alpha_2}} (...
{\rm Coeff}_{x_n^{\alpha_n}}  \nn
&& \frac{g(x_i,e^{y_j})_{i,j=1}^n}
{\displaystyle{\prod_{i<j}}(e^{y_i}x_i-e^{y_j}x_j)^{p_{i,j}}
(e^{y_i}x_i-x_j)^{r_{i,j}}(x_i-e^{y_j}x_j)^{s_{i,j}}
(x_i-x_j)^{u_{i,j}}}))=\nn
&&={\rm Coeff}_{x_1^{\alpha_1}}({\rm Coeff}_{x_2^{\alpha_2}}(...
{\rm Coeff}_{x_{n-1}^{\alpha_{n-1}}}  \nn
&& \frac{h(x_i,e^{y_j})_{i,j=1}^{n-1}}{
\displaystyle{\prod_{i<j}}(e^{y_i}x_i-e^{y_j}x_j)^{p_{i,j}}
(e^{y_i}x_i-x_j)^{r_{i,j}}(x_i-e^{y_j}x_j)^{s_{i,j}}
(x_i-x_j)^{u_{i,j}}} )),
\eea
where $h$ is some Laurent polynomial. By induction
it follows that (\ref{induction}) is an element of $\mathbb{C}[[y_1,...,y_n]]$.
Thus $f_{\alpha}(e^{y_i})_{i=1}^n \in \mathbb{C}((y_1,...,y_n))$.
\epfv

\begin{definition}
{\em Let $u_i,v_j \in V$, for $i=1,...,n+1$, $j=1,...,n$ and $u'_{n+1} \in V'$.
We define an iterated $2n$--point function as a formal matrix coefficient:
\begin{equation} \label{corre0}
\langle u'_{n+1}, X(Y[u_1,y_1]v_1,x_1)\cdots X(Y[u_n,y_n]v_n,x_n)u_{n+1}
\rangle,
\end{equation}
where
\begin{equation} \label{zhuv}
Y[u,y]=e^{{\rm wt}(u)y}Y(u,e^{y}-1).
\end{equation}
}
\end{definition}

Vertex operators of the form (\ref{zhuv}) have been introduced by Y. Zhu in
\cite{Zh1}. A quadruple $(V,Y[ \  , \  ], {\bf 1}, \tilde{\omega})$,
where $\omega=\omega-\frac{c}{24}$,
is a vertex operator algebra isomorphic to $(V, Y(\ ,\ ), {\bf 1}, \omega)$.
Operators of the type $X(Y[u,y]v,x)$ were introduced in \cite{Le1}, \cite{Le2}
in connection with the Riemann $\zeta$--function (cf. \cite{M}).

\subsection{Associativity for iterated $2n$--point functions}

The aim is to relate (\ref{corre0}) and (\ref{genform1}).
Let us consider the simplest case when $n=1$.
We know (cf. \cite{M} ) that
$$X(Y[u,y]v,x)=\nb X(u,e^y x)X(v,x) \nb + X(Y^+[u,y]v,x),$$
where $$Y^+[u,y]=\sum_{n \geq 0} u(n)(e^y-1)^{-n-1}$$
and the normal ordering  $\nb \ \nb$ for $X$--operators
is the one introduced in Section \ref{trivial}.
Hence
\bea \label{screw1}
X(Y[u,y]v,x_1)& = &\nb X(u,e^y x)X(v,x) \nb + X(Y^+[u,y]v,x) \nn
&=& \nb X(u,e^y x)X(v,x) \nb+\sum_{i \geq 0} (e^y)^{wt(u)}\frac{
X(u(i)v,x)}{(e^y-1)^{i+1}}.
\eea


On the other hand
\bea \label{screw2}
&& X(u,t x) X(v,x)=\nn
&& = \nb X(u,tx) X(v,x) \nb + [Y^+((tx)^{L(0)} u,t x),X(v,x)]
\eea


Now we work out the second term in (\ref{screw2}). \bea \label{screw4} &&
[Y^+((tx)^{L(0)}u,tx),X(v,x)]= \nn && = {(tx)}^{{\rm wt}(u)}
x^{{\rm wt}(v)} \lim_{x_1 \to x} {\rm Sing}_{x_1} \ [Y^+(u,tx_1),
Y(v,x)]= \nn && = {(tx)}^{{\rm wt}(u)} x^{{\rm wt}(v)} \lim_{x_1
\to x} {\rm Sing}_{x_1} {\rm Res}_{x_0} x^{-1} \delta
\left(\frac{tx_1-x_0}{x} \right) Y(Y(u,x_0)v,x) = \nn &&
={(tx)}^{{\rm wt}(u)} x^{{\rm wt}(v)} \lim_{x_1 \to x} {\rm
Sing}_{x_1} \sum_{i \geq 0} x^{-1} e^{-x_0
\frac{\partial}{\partial(tx_1)}} \delta \left(\frac{tx_1}{x}
\right) Y(u(i)v,x)=\nn && = {(tx)}^{{\rm wt}(u)} x^{{\rm wt}(v)}
\lim_{x_1 \to x} \sum_{i \geq 0} \frac{1}{(tx_1-x)^{i+1}}
Y(u(i)v,x)=\nn && ={t}^{{\rm wt}(u)} x^{{\rm wt}(u)+{\rm wt}(v)}
\sum_{i \geq 0} \frac{x^{-i-1}}{(t-1)^{i+1}} Y(u(i)v,x)= \nn && =
\sum_{i \geq 0} {t}^{{\rm wt}(u)} x^{{\rm wt}(u)+{\rm wt}(v)-i-1}
\frac{1}{(t-1)^{i+1}} Y(u(i)v,x)= \nn && = \sum_{i \geq 0}
\frac{{t}^{{\rm wt}(u)}}{(t-1)^{i+1}} X(u(i)v,x), \eea
where ${\rm Sing}_{x_1}$ stands for the regular part with respect
to $x_1$ and $\displaystyle{\lim_{x_1 \to x}}$ stands for the
formal substitution $x_1 \mapsto x$ (cf. \cite{FLM}).

It is obvious also that \bea \label{ty} && \langle  u_2', \nb
X(u,tx)X(v,x) \nb u_2 \rangle \in \mathbb{C}[t,x,t^{-1},x^{-1}],
\nn && \langle  u_2', \nb X(u,e^y x)X(v,x) \nb u_2 \rangle \in
\mathbb{C}[e^y,x,e^{-y},x^{-1}], \eea for $u_2 \in M$ and $u_2'
\in M'$.
%


If we combine (\ref{screw1}), (\ref{screw2})
and (\ref{screw4}) we obtain
\begin{proposition} \label{genus01}
There exist a polynomial $g$ such that
$$ \langle  u'_2, X(u_1,t x) X(v_1,x) u_2 \rangle
=\frac{g(tx,t)}{t^r x^s (t-1)^p},$$
for some $p,r,s \in \mathbb{N}$, $g(x,y) \in \mathbb{C}[x,y]$.
Moreover,
$$ \langle  u'_2, X(Y[u_1,y]v_1,x) u_2 \rangle=
\frac{g(e^yx,x)}{(e^y)^r x^s (e^y-1)^p}.$$
In particular, when $x$ and $y$ are specialized to
be complex variables  $$ \langle  u'_2, X(Y[u_1,y]v_1,x) u_2 \rangle $$
converges to a holomorphic function inside
$|x|>0$, $0<|y|<2 \pi $.
\end{proposition}
{\em Proof:}
We may assume that all vectors are homogeneous.
Relations (\ref{ty}), (\ref{screw2}) and (\ref{screw4}) imply
that $$\langle  u'_2, X(Y[u_1,y]v_1,x) u_2 \rangle=
\frac{g(e^yx,x)}{(e^y)^r x^s (e^y-1)^p}.$$
The second statement follows from the above discussion
and the fact that $\frac{1}{e^y-1}$, as a complex
function, has a Laurent expansion $\frac{1}{y}+ \cdots$
which converges inside $0<|y|<2 \pi$.
This expansion coincides with the expansion
$$\frac{1}{e^y-1}=
\frac{1}{y(1+\frac{e^y-1-y}{y} ) }=\frac{1}{y}-\frac{e^y-1-y}{y}+\cdots $$
exhibited according to our conventions.


\epfv


Let us consider the general case ($n \geq 2$).
For every $n \in \mathbb{N}$ we define domains
\bea
&& \Omega_{n,n}=\{ (x_1,...,x_n,t_1,...,t_n) \in \mathbb{C}^{2n}: \ \nn
&& |t_1x_1|>|x_1|>|t_2x_2|>|x_2|>....>|t_{n-1}x_{n-1}|>|x_n|>0 \}.
\eea
It is easy to see that (\ref{genform}) is
analytic inside $\Omega_{n,n}$.
Note that  $|t_i|>1$, holds for every $i=1,...,n$ and
$(x,t):=(x_1,...,x_n,t_1,...,t_n) \in \Omega_{n,n}$.
We have the following theorem (that we call {\em associativity})

\begin{theorem} \label{main1}
There exist a formal series $f(x_i,t_j)_{i,j=1}^n$ of the form
(\ref{genform}) such that
\begin{equation} \label{1a}
\langle u'_{n+1}, X(u_1,x_1 t_1)X(v_1,x_1) \cdots X(u_n,x_nt_n)
X(v_n,x_n) u_{n+1} \rangle=f(x_i,t_j)_{i,j=1}^n.
\end{equation}
and
\begin{equation} \label{1b}
\langle u'_{n+1} X(Y[u_1,y_1]v_1,x_1) \cdots
X(Y[u_n,y_n]v_n,x_n)u_{n+1} \rangle =f(x_i,e^{y_j})_{i,j=1}^n,
\end{equation}
where $f(x_i,e^{y_j})_{i,j=1}^n$ is a series viewed according
to our expansion conventions.
In particular, as complex functions, (\ref{1a}) and (\ref{1b})
have the same meromorphic continuation onto $\mathbb{C}^{2n}$ if
we replace $e^{y_i}$ by $t_i$.
\end{theorem}
{\em Proof:}
We know that
\bea
&& \langle u'_{n+1}, X(u_1,x_1 t_1)X(v_1,x_1) \cdots X(u_n,x_nt_n)
X(v_n,x_n) u_{n+1} \rangle = \\
&& \langle  u'_{n+1}, (\nb X(u_1,x_1 t_1)X(v_1,x_1) \nb
+F_1(t_1,x_1)) \cdots \nn && \cdots (\nb X(u_1,x_1 t_1)X(v_1,x_1)
\nb +F_n(t_n,x_n)) u_{n+1} \rangle, \nonumber \eea where
$$F_i(t_i,x_n)= [Y^+((t_ix_i)^{L(0)}u,tx),X(v,x_i)]$$
On the other hand
\bea \label{screw5}
&& \langle
u'_{n+1}, X(Y[u_1,y_1]v_1,x_1) \cdots X(Y[u_n,y_n]v_n,x_n)
u_{n+1} \rangle= \\
&& = \langle u'_{n+1}, ( \nb X(u_1,e^{y_1}x_1) X(v_1,x_1) \nb +
f_1(e^{y_1},x_1) ) \ldots  \nn && \ldots ( \nb X(u_n,e^{y_n}x_n)
X(v_n,x_n) \nb + f_n(e^{y_n},x_n) ) u_{n+1} \rangle, \nonumber
\eea
where
$$f_i(e^{y_i},x_i)=X(Y^+[u_i,y_i]v_i,x_i),$$
and the (operator valued) series on the
right hand side are obtained according to
our expansion convention.

Because (\ref{1a}) is a sum of several terms of the form
%
\begin{equation} \label{morgenpole}
\langle u'_{n+1}, \nb X(u_1,x_1 t_1)X(v_1,x_1)
\nb \cdots \frac{{t_j}^{{\rm wt}(w)} X(w,x_j)}{(t_j-1)^m} \cdots \nb X(u_n,x_n t_n)X(v_n,x_n) \nb u_{n+1} \rangle,
\end{equation}
for some $m \geq 0$,
the only thing we have to prove is that

\begin{equation} \label{morgenpole1}
\langle u'_{n+1}, \nb X(u_1,x_1 t_1)X(v_1,x_1)
\nb \cdots {X(w,x_j)} \cdots \nb X(u_n,x_n t_n)X(v_n,x_n) \nb u_{n+1} \rangle,
\end{equation}
is an expansion of a meromorphic function with (possible) poles at
$x_i=0$, $t_i=0$, $x_i=x_j$, $t_ix_i=x_j$, $i \neq j$ and
$t_ix_i=t_jx_j$. If this is the case then (\ref{morgenpole1}) does
not contribute with additional poles at $t_i=1$ so we are allowed
to make a substitution $t_i=e^{y_i}$ in (\ref{morgenpole1}) and
use the same expansion conventions (which are applicable !) and
obtain
\begin{equation} \label{morgenpole2}
\langle u'_{n+1}, \nb X(u_1,x_1 e^{y_1})X(v_1,x_1)
\nb \cdots {X(w,x_j)} \cdots \nb X(u_n,x_n e^{y_n})X(v_n,x_n) \nb u_{n+1} \rangle.
\end{equation}
But we know that (\ref{screw5}) is a sum of several terms like
\begin{equation} \label{morgenpole3}
\langle u'_{n+1}, \nb X(u_1,x_1 e^{y_1})X(v_1,x_1)
\nb \cdots \frac{e^{y_j {\rm wt}(w) } X(w,x_j) }{(e^{y_j}-1)^m} \cdots \nb X(u_n,x_n e^{y_n})X(v_n,x_n) \nb u_{n+1} \rangle,
\end{equation}
so we have the proof. \\
\noindent {\em Claim:} The expression (\ref{morgenpole}) is an expansion of a rational
function with no poles at $t_i=1$. \\
\noindent {\em Proof:}
If we use
$$ \nb X(u_i,x_i t_i)X(v_i,x_i) \nb=X^-(u_i,x_i t_i)X(v_i,x_i)+
X(v_i,x_i)X^+(u_i,x_i t_i),$$
for $i=1$ we get
\bea \label{screw4a}
&& \langle  u'_{n+1}, \nb X(u_1,x_1 t_1)
X(v_1,x_1) \nb \cdots X(w,x_j) \cdots \nb X(u_n,x_n t_n)X(v_n,x_n) \nb u_{n+1}
\rangle= \nn && = \langle  u'_{n+1}, X^+(u_1,x_1 t_1)X(v_1,x_1)\nb
\cdots X(w,x_j) \cdots \nb X(u_n,x_n t_n)X(v_n,x_n) \nb u_{n+1} \rangle+ \nn
&& \langle  u'_{n+1}, X(v_1,x_1)X^-(u_1,x_1 t_1)\cdots \nb X(u_n,x_n
t_n)X(v_n,x_n) \nb u_{n+1} \rangle.
\eea
After we apply the same
formula, for every $i$, we obtain several summons. A generic term
is
\begin{equation} \label{screw6}
\langle  u'_{n+1}, X(v_1,x_1)X^-(u_1,x_1 t_1) \cdots X(w_j,x_j) \cdots X^+(u_n,x_n t_n)X(v_n,x_n)
u_{n+1} \rangle.
\end{equation}
It is important to notice that $X^+(u_i,t_ix_i)$ is always to the right
of $X(v_i,x_i)$ and $X^-(u_i,t_ix_i)$ is to the left of $X(v_i,x_i)$.

%
Now, by using the commutator formula we can move $X^+(u_i,t_ix_i)$ all the way to the
right (such that it acts on $w_{n+1}$) and $X^-(u_i,x_i)$ all the way to the left.
This procedure will produce poles only at wanted places, i.e., no new poles
at $t_i=1$. After we apply  $X^+(u_i,t_ix_i)$ to $w_{n+1}$ we obtain only finitely many terms.
The same thing happens with $X^-(u_i,x_i)$. We can repeat the same procedure
for remaining terms. This is a finite procedure and at the end we end up with a rational function
with no poles at $t_i=1$.


The proof follows.

\epfv

We define an {\em $n$--point function} to be a formal matrix
coefficient: \bea \label{screw9} \langle u'_{n+1},
o(Y[u_1,y_1]v_1)\cdots o(Y[u_n,y_n]v_n)u_{n+1} \rangle, \eea where
$o(a)=a({\rm wt}(a)-1)$, i.e., we extract the zeroth coefficient
in (\ref{corre0}) with respect to the expansion in terms of
$x_i$'s.  In parallel with the previous construction we will
consider: \bea \label{screw10} \langle u'_{n+1},
o(X(u_1,t_1x_1)X(v_1,x_1)) \cdots o(X(u_n,t_nx_n)X(v_n,x_n))
u_{n+1} \rangle, \eea where
$$o(X(u,t x)X(v,x)):={\rm Res}_{x}x^{-1} X(u,t x)X(v,x).$$

\noindent For $n \geq 1$, let
$$\Omega_{n,0}=\{ (t_1,...,t_n) \in \mathbb{C}^n : |t_i|>1,
i=1,...,n \}.$$
Then we have
\begin{corollary}
The expression
\bea \label{screw11}
&& f(t_1,...,t_n)= \\
&& = \langle u'_{n+1}, o(X(u_1,t_1x_1)X(v_1,x_1))\cdots
o(X(u_n,t_nx_n)X(v_n,x_n)) u_{n+1} \rangle. \nonumber
\eea
converges to a holomorphic function
inside  $\Omega_{n,0}$ and it has a meromorphic extension
to $\mathbb{C}^n$ with poles at $t_i=1$ and $t_i=0$.
\bea \label{screw12}
g(y_1,...,y_n)=\langle u'_{n+1}, o(Y[u_1,y_1]v_1) \cdots
o(Y[u_n,y_n]v_n) u_{n+1} \rangle
\eea
converges inside $0<|y_i|<2 \pi$. Moreover, as meromorphic functions,
\begin{equation}
g(y_1,...,y_n)=f(t_1,...,t_n)|_{t_1=e^{y_1},...,t_n=e^{y_n}},
\end{equation}
\end{corollary}
{\em Proof:}
From the Proposition \ref{welld} we know that (\ref{screw12})
is well--defined.
Here we make the following additional observation:
$${\rm Coeff}_{x_1^0 \cdots x_n^0} \left \{
\frac{g(x_i,t_j)_{i,j=1}^n}{\displaystyle{\prod_{i <j}}
(t_ix_i-x_j)^{r_{i,j}} (x_i-t_jx_j)^{s_{i,j}}(x_i-x_j)^{u_{i,j}}} \right\}
$$
is a Laurent polynomial in $t_i$, $i=1,...,n$.
Because of (\ref{genform1}), it follows that (\ref{screw11})
converges inside $|t_i|>1$.
Similarly,
$${\rm Coeff}_{x_1^0 \cdots x_n^0} \left \{
\frac{g(x_i,e^{y_j})_{i,j=1}^n}{\displaystyle{\prod_{i <j}}
(e^{y_i}x_i-x_j)^{r_{i,j}} (x_i-e^{y_j} x_j)^{s_{i,j}}(x_i-x_j)^{u_{i,j}}}
\right \}
$$
is (the same) Laurent polynomial in $t_i=e^{y_i}$, $i=1,...n$ and it
converges inside $0<|y_i|<2\pi$.
Now, Theorem \ref{main1} implies the last statement in the corollary.

\epfv


\begin{remark}
{\em In Proposition \ref{genus01} we did not
use the associativity (as stated in
\cite{FHL}). One might think that there is an easier proof (cf. \cite{Le1}--\cite{Le3}):
\bea \label{nonrig}
&& X(u,t x)X(v,x)=(t x)^{{\rm wt}(u)}x^{{\rm wt}(v)}Y(u,t x)Y(v,x)
\sim \nn
&& (t_1x)^{{\rm wt}(u)}x^{{\rm wt}(v)}Y(Y(u,tx-x)v,x) \sim \nn
&&(t_1x)^{{\rm wt}(u)}  X(x^{-L(0)}Y(u,tx-x)x^{L(0)} v,x) \sim \nn
&&X(t^{L(0)}Y(u,t-1)v,x).
\eea
Even though (\ref{nonrig}) implies the ``right'' answer (after we
substitute $e^{y}$ for $t$), the proof
is not rigorous ! Note that the term $X(t^{L(0)}Y(u,t-1)v,x)$
is a non-rigorous expression (this can be seen if we extract
zeroth coefficient of $t$) compared with the left hand side which is
well-defined. However, if we use matrix coefficients this construction
can be made completely rigorous. }
\end{remark}

\subsection{$q$--traces for the iterated $2n$--point functions}

In the previous section we showed that
two correlation functions are related by a simple substitution (Theorem \ref{main1}).
Here we extend the same result for the corresponding $q$--traces.
Let us consider (cf. \cite{Zh1}) a  formal series ($q$--trace):
\bea \label{zhugen}
{\rm tr}|_{M} X(a_1,x_1) \cdots X(a_n,x_n) q^{L(0)} \in
q^{h}\mathbb{C}[[x_1^{\pm 1},...,x_n^{\pm 1},q]],
\eea
$a_i \in V$, $i=1,...,n$, where $M$ is a $V$--module which admits
$\mathbb{N}$--grading,
i.e., there is $h \in {\mathbb C}$ such that ${\rm Spec}(L(0)) \in
h+\mathbb{Z}_{\geq 0}$.
For $u_i,v_i \in V$, $i=1,...,n$
we define a $q$--trace for an iterated $2n$--point function as
a formal expression
\begin{equation} \label{2npoint}
{\rm tr}|_{M} X(Y[u_1,y_1]v_1,x_1) \cdots X(Y[u_n,y_n]v_n,x_n) q^{L(0)}.
\end{equation}
Later $q=e^{2 \pi i \tau}$, where $\tau \in \mathbb{H}$ ($\mathbb{H}$
is the upper half-plane).

It is a novelty, comparing to (\ref{zhugen}), to consider correlation
functions of the form (\ref{2npoint}).
Now the vector $a_i$ is replaced by the formal expression $Y[u_i,y]v_i$. Thus one can think
of (\ref{2npoint}) as a $q$--trace
attached to a certain completion of $V$.

We will work intensively with another related $q$--trace:
\bea \label{2nxs}
{\rm tr}|_{M} X(u_1,t_1x_1)X(v_1,x_1) \cdots X(u_n, t_nx_n)X(u_n, x_n)
q^{L(0)}.
\eea

By the definition:
\bea \label{trxs}
&& {\rm tr}|_M
X(u_1,t_1x_1)X(v_1,x_1) \cdots X(u_n,t_nx_n)X(v_n,x_n) q^{L(0)}=
\eea
$$ \sum_{i \in I }
\langle
w_{i}',X(u_1,t_1x_1)X(v_1,x_1) \cdots X(u_n,t_nx_n)X(v_n,x_n)w_{i}
\rangle
q^{\rm{wt}(w_{i})},$$
where $\{w_i \}_{i \in I}$ is some (or any) homogenuous basis for $M$ such that
$\langle w'_{i},w_{i} \rangle =1$.
We know that
$$ \langle w_{i}',X(u_1,t_1x_1)X(v_1,x_1) \cdots X(u_n,t_nx_n)X(v_n,x_n)w_{i}
\rangle $$ converges to an analytic function in the domain
$$|t_1x_1|>|x_1|>...>|t_nx_n|>|x_n|>0,$$
but {\em a priori}  we do not know
whether (\ref{2nxs}) has a domain of convergence.
To prove the convergence it suffices
(according to  \cite{Zh1})  to check  that for every $v \in V$ the formal 1--point function
$${\rm tr}|_M X(v,x) q^{L(0)}$$
converges  to a holomorphic function inside $|q|<1$.
This statement is very strong and it depends on the internal structure of
the vertex operator algebra ($C_2$--condition for instance).
If this is the case then one can show that (\ref{trxs}) converges to
a holomorphic function
inside
\begin{equation} \label{natdom}
1>|t_1x_1|>|x_1|> \ldots >|t_nx_n|>|x_n|>|q|,
\end{equation}
and it has a meromorphic extension to $(\mathbb{C}^\times)^{2n}$.

In our approach $1$--point function does not play a prominent role
(actually we start with $2$--point functions). Also we do
not assume rationality, $C_2$--condition etc..
Therefore we need to develop a right recursion procedure
for expressing  $q$--traces of iterated $2n$---point functions
in terms of $q$--traces of certain $(2n-2)$--point functions.
Then, by using this result we prove the convergence of the
iterated $2n$--point. We shall adopt this
approach at the very end.

From now on we will assume that (\ref{condition}) holds for all
pairs of vectors $u_i,v_j \in V$, $i=1,...,n$.
The following fact follows immediately from Theorem \ref{main1}.

\begin{corollary} \label{main2}
Let us denote (\ref{2nxs}) by $f(t,x,q)$ and (\ref{2npoint}) by
$g(y,x,q)$. Then \bea \label{fg} &&g(y,x,q)=f(e^y,x,q). \eea If we
assume that $x,y$ and $t$ are complex variables and $q=e^{2 \pi i
\tau}$ then (\ref{fg}) holds, provided that both sides are
absolutely convergent in a certain domain.
\end{corollary}

The previous corollary is not so useful since we do not have any
information about poles.

In the case $n=1$ we have the following description \cite{Zh1} (for the notation
and a different proof see Appendix A).
\begin{theorem} \label{th4}
Formally,
\bea \label{c2zhu}
&& {\rm tr}|_M X(Y[u,y]v,x)=\\
&& \sum_{m \geq 1}\bar{\wp}_{m+1}(y,q) {\rm tr}|_M X(u[m]v,x) q^{L(0)}+
{\rm tr}|_M o(u)o(v)q^{L(0)}, \nonumber
\eea
where $\bar{\wp}_{m+1}(y)$ (defined in the Appendix B)
are considered as elements from
$\mathbb{C}((y))[[q]]$.
\end{theorem}

\begin{remark} \label{multival}
{\em All results obtained in the previous section carry out
(in a straightforward manner) for vertex operator
superalgebras (cf. \cite{KW}).
Because of the $\frac{1}{2}\mathbb{Z}$--grading,
correlation functions are multivalued functions and therefore
the statements in Theorem \ref{main1} and Corollary
\ref{main2} are ambiguous.
If we assume that $p(u_i)=p(v_i)$ (where $\epsilon$ is
the sign), in (\ref{2nxs}) possible multivaluedness stems
only from the term
$$\sqrt{t_{i_1}...t_{i_k}},$$
for some $1 \leq i_i <...< i_k \leq n$.
}
\end{remark}

\renewcommand{\theequation}{\thesection.\arabic{equation}}
\setcounter{equation}{0}

\section{$q$--traces and $\hat{\mathcal{D}}$; charged fermions}

\subsection{Charged free fermions}
In this section we exploit our general consideration
in the case of the vertex operator superalgebra stemming
from the pair of charged fermions (\cite{KP}, \cite{K}). Also, we discuss
an equivalent bosonic construction.

First we introduce free (charged) fermions.
We consider a Lie superalgebra generated by
$\psi_n$ and $\psi^*_m$, $n, m \in {\mathbb Z}+1/2$, and $1$
with commutation relations
$$[\psi_n,\psi^*_m]=\delta_{m+n,0}. \ \  [\psi^*_n,\psi^*_m]=[\psi_n,\psi_m]=0.$$
We denote the corresponding  Fock space
by ${\cal F}$ (\cite{K}) spanned by elements
\begin{equation} \label{span}
\psi^*(-n_k) \cdots \psi^*(-n_1)\psi(-m_l) \cdots \psi(-m_1){\bf 1},
\end{equation}
where $n_i ,m_j \in {\mathbb N}+\frac{1}{2}$ and
$n_k> \ldots >n_1\geq 1/2$, $m_l> \ldots >m_1 \geq \frac{1}{2}$.

Let us define the following vertex operators (fermionic fields):
$$Y(\psi(-1/2){\bf 1},x)=\sum_{n \in \mathbb Z} \psi(n+1/2)x^{-n-1}$$
and
$$Y(\psi^*(-1/2){\bf 1},x)=\sum_{n \in \mathbb Z} \psi^*(n+1/2)x^{-n-1}.$$
We write shorthand $\psi=\psi(-1/2){\bf 1}$ and $\psi^*=\psi^*(-1/2){\bf 1}$.
Since the fermionic fields are local and generate the space ${\cal F}$
we can equip ${\cal F}$ with a structure of vertex operator superalgebra
(see \cite{FKRW}) with the conformal vector
$$\omega=\frac{1}{2}\psi^*(-3/2)\psi(-1/2){\bf
1}+\frac{1}{2}\psi(-3/2)\psi^*(-1/2){\bf 1}.$$
This vertex operator superalgebra has the {\em charge}
decomposition
$${\cal{F}}=\bigoplus_{n \in \mathbb{Z}} {\cal{F}}_n,$$
with respect to the operator
$$o(\nb X(\psi,x)X(\psi^*,x) \nb).$$
Then ${\cal{F}}_0$ is a vertex operator algebra and
every ${\cal{F}}_n$ is an ${\cal{F}}_0$--module.

The projective representation of the Lie
algebra $\mathcal{D}$ of differential
operators on $\mathbb{C}^\times$ (see \cite{KR1}, \cite{FKRW})
, studied in \cite{BO}, can be interpreted
in the language of the vertex operator
algebras in the following way.

Consider the following vectors:
\begin{equation} \label{elements}
\psi(-i-1/2)\psi^*(-j-1/2){\bf 1},
\end{equation}
$i, j \in {\mathbb N}$.
Then ${\cal D}$ has a projective representation
in terms of the Fourier coefficients of the vectors (\ref{elements}) (see \cite{FKRW}).
Let ${\cal D}_0$ be the Cartan subalgebra of ${\cal
D}$ (see \cite{KR1}) spanned by $L^{(r)}_0=(t \frac{d}{dt})^r$, $r {\geq 0}$.

We define an action of ${\cal D}_0$ in terms
of generating functions (we choose a particular lifting) in the
following way:
$$L^{(r)}_0 \mapsto \bar{L}^{(r)}(0),$$
where  $\bar{L}^{(r)}(0)$ is the coefficient of
$x^0y^r$ in
\begin{equation} \label{dfermgen}
\nb X(\psi, e^y x) X(\psi^*,x) \nb +\frac{e^{y/2}}{e^y-1}.
\end{equation}

According to \cite{Bl}, \cite{Le1} and \cite{M} (\ref{dfermgen})
corresponds to a new normal ordering.
The operator $\bar{L}^{(r)}(0)$ introduced
in  \cite{BO} coincide with $\tilde{D}_r$ for every $r \in {\mathbb N}$.
It is easy to see that operators $\bar{L}^{(r)}(0)$
act semisimply on ${\cal F}$ and therefore we can
define a generalized character of ${\cal F}$ (cf. \cite{FKRW}):
\begin{equation}
{\rm tr}|_{\cal{F}} \prod_{i=0}^\infty q_i^{\bar{L}^{(i)}(0)},
\end{equation}
where $q_i=e^{2 \pi  i \tau_i}$, for $i \geq 1$ where
$\tau_1 \in \mathbb{H}$, and $q_0$, $\tau_i$, $i \geq 2$
are formal variables.
Notice that $\bar{L}^{(0)}(0)$ is not associated to any
particular element of $\hat{\mathcal{D}}$.

It is important to notice that (\ref{dfermgen}) can be
written as
\begin{equation} \label{dferm}
X(Y[\psi,y]\psi^*, x).
\end{equation}

\begin{remark}
{\rm One can derive a closed formula for the
commutator
\begin{equation}
[X(Y[\psi,y_1]\psi,x_1), X(Y[\psi,y_2]\psi^*,x_2)],
\end{equation}
by using the Jacobi identity from \cite{M} or (\ref{dfermgen}).
We leave this (non-trivial) exercise to the reader.}
\end{remark}

Then
\begin{eqnarray} \label{bulk}
&&{\rm tr}|_{\cal{F}} \prod_{i=0}^\infty q_i^{\bar{L}^{(i)}(0)}=
\prod_{i=1}^{\infty} q_{2i-1}^{\zeta(-2i+1,1/2)} \nn
&& \prod_{n \geq 0}
(1+q_0q_1^{n+1/2}q_2^{(n+1/2)^2} \cdots
)(1+q_0^{-1}q_1^{n+1/2}q_2^{-(n+1/2)^2} \cdots ).
\end{eqnarray}
This corresponds to the formula (1.12) in \cite{BO}, where the same
generating function is denoted by $\Omega(\tau_0,\tau_1,...)$.
The infinite product has an expansion as a Laurent series of
$q_0$. We  denote by $V(\tau_1, \tau_2...)$
the coefficient of $q_0^0$ of the generating function (\ref{bulk}).
$\Omega$ and $V$  are {\em quasi-modular forms}
(for the definition of quasi--modularity see \cite{BO}, \cite{M}) of the
weight $0$ and $-\frac{1}{2}$ respectively.

\subsection{Boson-fermion correspondence} \label{bsc}

Let $L$ be a rank one lattice with a generator $\alpha$ such that
$<\alpha,\alpha>=1$.
We denote by $\hat{L}$ the non--trivial central extension
of $L$
$$1 \rightarrow \{\pm 1 \} \rightarrow \hat{L}
\stackrel{-}{\rightarrow} L \rightarrow 0, $$
such that $ab=(-1)^{<\bar{a},\bar{b}>}ba$, $a,b \in \hat{L}$.

As in  \cite{DL} (or \cite{K}) we equip the space
\begin{equation} \label{dong}
V_L \cong M(1) \otimes {\mathbb C}[L],
\end{equation}
with a structure of vertex operator superalgebra.
Then
\begin{equation} \label{lattice}
X(a^{\pm 1},x)={\rm exp} \left(\sum_{n >0}\frac{\alpha(-n)}{-n}x^{n}
\right) {\rm exp} \left(\sum_{n >0}\frac{\alpha(n)}{n}x^{-n}\right)a^{\pm 1} z^{\pm \bar{a}},
\end{equation}
where $\bar{a}=\alpha$.
The boson-fermion
correspondence  (cf. \cite{K}) is a vertex operator
superalgebra isomorphism
\begin{equation}
\Psi : \ {\cal F} \  \rightarrow \ V_L,
\end{equation}
with the following properties: \bea && \psi(-n+1/2) \cdots
\psi(-1/2){\bf 1} \mapsto e^{n\alpha} \nn && \psi^*(-n+1/2) \cdots
\psi^*(-1/2){\bf 1} \mapsto e^{-n\alpha}, \nn && \omega \mapsto
\frac{1}{2}\alpha(-1)^2{\bf 1}. \eea In particular,
$\psi(-1/2){\bf 1} \mapsto e^{\alpha}$ and $\psi^*(-1/2){\bf 1}
\mapsto e^{-\alpha}$.

Now (\ref{dferm}) can be expressed in the following form:
\begin{equation} \label{sugawara}
X(Y[\psi, y]\psi^*, x)=\Psi^{-1} X(Y[e^{\alpha}, y]
e^{-\alpha}, x) \Psi,
\end{equation}
since
$$\Psi \left( Y[\psi,y]\psi^* \right)=Y[e^\alpha,y]e^{-\alpha},$$
which is a consequence of Proposition 4.1 in \cite{M}.

If we denote by $\tilde{L}^{(r)}(0)$ the coefficient of
$x^0y^r$ in
$$X(Y[a, y]a^{-1}, x),$$
then
\begin{equation}
{\rm tr}|_{\cal{F}} \prod_{i=0}^\infty q_i^{\bar{L}^{(i)}(0)}=
{\rm tr}|_{V_L} \prod_{i=0}^\infty q_i^{\tilde{L}^{(i)}(0)}.
\end{equation}
In particular if we consider only the coefficient of $q_0^0$ , i.e.
$V(\tau_1,\tau_2,...)$, we get
\begin{equation}
{\rm tr}|_{{\mathcal{F}}_0} \prod_{i=0}^\infty q_i^{\bar{L}^{(i)}(0)}
= {\rm tr}|_{M(1)} \prod_{i=0}^\infty q_i^{\tilde{L}^{(i)}(0)},
\end{equation}
where $M(1)$ is identified with $M(1) \otimes 1 \subset V_L$ by
means of (\ref{dong}).

\begin{remark}
{\em Generalized characters discussed in this section
already appeared in the literature.
In \cite{Di}, \cite{Di1} similar generating functions related
to mirror symmetry on the torus resemble (\ref{bulk}).
Also they appear (with finitely many $\bar{L}^{(r)}(0)$'s) in  \cite{KZ}.}
\end{remark}

\subsection{Bloch-Okounkov's $n$--point function}

In \cite{BO} Bloch and Okounkov introduced  another type of generating
functions. Let us associate to the generalized
character $V(\tau_1,\tau_2,...)$ an {\em $n$--point function}
in the following way. First we define operators
$$\sigma(y)=\frac{1}{y}+\frac{1}{2 \pi i}
\sum_{r=1}^{\infty} \frac{\partial}{\partial \tau_r} \frac{y^r}{r!}.$$
Then the (formal) $n$-point function is defined as
\begin{equation}
{\cal F}(y_1,...,y_n)=\eta(q_1) \sigma(y_1) \cdots \sigma(y_n) V(\tau_1,
\tau_2,...)|_{\tau_2=\tau_3=...0}.
\end{equation}
We mention that $u_i$'s is a formal variables
in contrast to \cite{BO} where $u_i$ is a
complex variable.

\begin{proposition} \label{deltax}
We have
\begin{eqnarray}
&& \eta(q_1) \sigma(y_1) \cdots \sigma(y_n) V(\tau_1,
\tau_2,...)|_{\tau_2=\tau_3=...0}= \nn
&& {\rm Coeff}_{x^0}
{\rm tr}}|_{{\cal F}_0}{ \eta(q_1) X(Y[\psi,y_1]\psi^*, x_1) \cdots X(Y[\psi ,y_n]\psi^*
,x_n)q_1^{\bar{L}(0)}=\nn
&&{\rm tr}|_{{\cal F}_0}  \eta(q_1) o(Y[\psi,y_1]\psi^*) \cdots o(Y[\psi ,y_n]\psi^*) q_1^{\bar{L}(0)},
\end{eqnarray}
where $\bar{L}(0)=L(0)-\frac{1}{24}$.
\end{proposition}
{\em Proof:}
The second equality is trivial since
$$X(u,x)=\sum_{n \in \mathbb{Z}} u({\rm wt}(u)-n-1) x^{-n}.$$
For the first equality we prove only in the case $n=1$.
Let us recall
$$V(\tau_1,\tau_2,...)={\rm tr}|_{\mathcal{F}_0} \prod_{i=1}^{\infty}
q_i^{\bar{L}^{(i)}(0)}.$$ 
Then 
\bea \label{1point} &&
\eta(q_1) \sigma(y)V(\tau_1, \tau_2,...)|_{\tau_2=\tau_3=...0}=\nn
&& {\rm tr}|_{\mathcal{F}_0} \ \eta(q_1)(
\frac{1}{y}+\sum_{r=1}^{\infty} {\bar L}^{(r)}(0)
\frac{y^{r}}{r!})q_1^{\bar{L}(0)}=\nn 
&& {\rm tr}|_{\mathcal{F}_0} \
\eta(q_1) \biggl(\frac{1}{y}-\sum_{r=1}^{\infty}
\zeta(-r,1/2)\frac{y^r}{r!}+o(\nb X(\psi, e^{y}x)X(\psi, x) \nb)
\nn && -o(\nb X(\psi, x)X(\psi, x) \nb) \biggr) q_1^{\bar{L}(0)}. \eea
Since the charge operator $o(\nb X(\psi, x)X(\psi, x) \nb)$ acts
as zero on ${\cal F}_0$ and
$$\frac{1}{y}-\sum_{r=1}^{\infty}
\zeta(-r,\frac{1}{2})\frac{y^r}{r!}=\frac{e^{y/2}}{e^y-1},$$
the proof follows from (\ref{dferm}).
\epfv

Now let us apply the results from Section 2.4  in the case of
vertex operator superalgebra $\mathcal{F}$ (with the obvious
super--modification). Fix $u_i=\psi$, $v_i=\psi^*$, $i=1,..,n$.
Notice that (unlike in Section 2.4.) we do not
consider an ${\cal F}$--module but rather a subalgebra of ${\cal F}$.
Since $X(\psi,t_ix_i)X(\psi^*,x_i)$ acts on ${\cal F}_0$ we can compute
the $q$--trace
\bea \label{2rxs}
{\rm tr}|_ {{\cal F}_0}X(\psi,t_1x_1)
X(\psi^*,x_1) \cdots X(\psi,t_nx_n)X(\psi^*,x_n) q^{\bar{L}(0)},
\eea
and a corresponding coefficient of $x^0$, multiplied by $\eta(q)$:
\bea \label{2nxs1}
&& G(t_1,...,t_n):= \\
&& = \eta(q) {\rm tr}|_{{\cal F}_0} o(X(\psi,t_1x_1)X(\psi^*,x_1))
\cdots o(X(\psi,t_nx_n)X(\psi^*,x_n)) q^{\bar{L}(0)}. \nonumber
\eea

\begin{proposition} \label{main3}
Suppose that (\ref{2nxs1}) converges, in some domain,
to a multivalued function $G(t_1,...,t_n)$.
Then, as complex functions,
$${\cal F}(y_1,...,y_n)=G(t_1,...,t_n)|_{t_i^{n/2}=e^{ny_i/2}}$$
for every $y$ provided that the right hand side is convergent.
\end{proposition}
{\em Proof:}
Since $$\psi(n)\psi^*=\delta_{n,1/2}{\bf 1},$$
for $n \geq \frac{1}{2}$ and
$$X(\psi,t_ix_i)X(\psi^*,x_i)= \nb X(\psi,t_i x_i)X(\psi^*,x_i) \nb +\frac{t^{1/2}}{t-1},$$
$$X(Y[\psi,y_i]\psi^*,x_i)= \nb X(\psi,e^{y_i} x_i)X(\psi^*,x_i) \nb +\frac{e^{y_i/2}}{e^{y_i}-1},$$
we can apply Theorem 2.1 and Corollary 2.1 (with minor modifications).
The correlation function $G$ converges to a multi valued function
(cf. Remark \ref{multival}). However, we can choose a branch
$t_i^{n/2}=e^{ny_i/2}$. Therefore the Proposition holds.
\epfv

We will need the following trivial (but important) observation.
\begin{lemma} \label{trivialemma}
Let $\varphi : V_1 \rightarrow V_2$ be a
vertex operator superalgebra isomorphism
and $W$ a subalgebra of $V_1$   stable under
$X(u_1,x_1) \cdots X(u_n,x_n)$. Then
\bea \label{traceiso}
&& {\rm tr}|_W X(u_1,x_1) \cdots X(u_n,x_n) q^{L_{V_1}(0)}=\nn
&& {\rm tr}|_{\varphi(W)} X(\varphi(u_1),x_1) \cdots X(\varphi(u_n),x_n)q^{L_{V_2}(0)}.
\eea
\end{lemma}
{\em Proof:} It follows immediately form the property of the trace
$${\rm tr}|_{\varphi(W)} \ \varphi X(u_1,x_1) \cdots X(u_n,x_n)
\varphi^{-1} q^{L_{V_2}(0)}=
{\rm tr}|_M \  X(u_1,x_1) \cdots X(u_n,x_n) q^{L_{V_1}(0)},$$
and the fact that
$$\varphi X(u_1,x_1) \varphi^{-1}=X(\varphi(u_1),x_1).$$
\epfv

\subsection{$q$--traces in the case of free fermions}

The boson-fermion correspondence implies
(cf. Lemma \ref{trivialemma}) that
\begin{equation} \label{mainpoint0}
{\rm tr}|_{{\cal F}}
X(\psi,t_1x_1)X(\psi^* ,x_1) \cdots X(\psi,t_n x_n)
X(\psi^*,x_n) q^{\bar{L}(0)},
\end{equation}
is equal to
\begin{equation} \label{mainpoint}
{\rm tr}|_{V_L}
X(a,t_1x_1)X(a^{-1},x_1) \cdots X(a,t_n x_n)
X(a^{-1},x_n) q^{\bar{L}(0)},
\end{equation}
and the corresponding trace for the zero--charge
subspace $\mathcal{F}_0$ is equal to (\ref{2rxs}).

Now we need some notation. Let $(q)_{\infty}=
\displaystyle{\prod_{i \geq 1}^\infty}(1-q^i)$
and
$$\theta(t) = \sum_{n \in \mathbb{Z}} q^{\frac{n^2}{2}} t^n \in
\mathbb{C}[[q^{1/2},t^{\pm 1}]].$$
If  $t=e^{2 \pi i v} \in \mathbb{C}^\times$
and $q=e^{2 \pi i \tau}$, $\tau \in \mathbb{H}$
then  $\theta(t)$ is the classical theta function.
Also we define another theta function in the infinite-product form
$$\theta_{11}(\tau,v)=\frac{\displaystyle{\prod_{i=1}^{\infty}(1-q^i
t)\prod_{i=1}^{\infty}}(1-q^i t^{-1})(t^{1/2}-t^{-1/2})}
{(q)^2_{\infty}}.$$
This is an entire function of $v$ for fixed $\tau \in \mathbb{H}$.
For a fixed $\tau$ we will write $\theta_{11}(v)$.
These theta functions can be thought as
functions of $t$ (instead of $v$)
but then these are multi valued function which we
will denote by $\theta(t)$ and $\theta_{11}(t)$.
The following transformation formulas hold:
$$\theta(q t)=q^{-1/2}t^{-1} \theta(t),$$
$$\theta_{11} (q t)=-q^{-1/2}t^{-1} \theta(t).$$

The following result has been known, in some form, for a while. 
For the completeness we include a proof here.
\begin{theorem} \label{mother}
\begin{enumerate}
\item[(a)]
\begin{equation} \label{zadnji}
{\rm tr}|_{\cal{F}}
X(\psi,t_1x_1)X(\psi^* ,x_1) \cdots X(\psi,t_nx_n)
X(\psi^* ,x_n) q^{L(0)}
\end{equation}
converges to a multi-valued holomorphic function in the domain
\begin{equation} \label{domain2}
|t_1x_1|>|x_1|>...>|t_nx_n|>|x_n|>|qt_1x_1|>0,
\end{equation}
\item[(b)]
\begin{equation} \label{ttq}
\eta(q) {\rm tr}|_{{\cal{F}}_0}
X(\psi ,t_1x_1)X(\psi^*,x_1) \cdots X(\psi,t_nx_n)
X(\psi^* ,x_n) q^{\bar{L}(0)}
\end{equation}
converges in the domain  (\ref{domain2}), and it has a
meromorphic extension to a double cover of $(\mathbb{C}^\times)^{2n}$ equals to
\bea \label{thetaproducts}
\frac{\displaystyle{\prod_{i <j}}
\theta_{11}(t_j x_j/t_ix_i) \theta_{11}(x_j/x_i)}
{\displaystyle{\prod_{i<j}} \theta_{11}(t_j x_j/x_i) \theta_{11}(x_j/t_ix_i)
\prod_{i=1}^n \theta_{11}(t_i)}. \ \
\eea

\end{enumerate}
\end{theorem}
{\em Proof:}
Let
$$G_n(t,x):=X(\psi, x_1t_1)X(\psi^*,t_1) \cdots X(\psi, x_nt_n)X(\psi^*,t_n).$$
First notice that,
because of the boson-fermion correspondence, we can replace \\
$X(\psi, x_1t_1)X(\psi^*,t_1)$ by $X(a, x_1t_1)X(a^{-1},t_1)$.
Also
\begin{equation} \label{prodtr}
{\rm tr}|_{M(1) \otimes {\mathbb{C}}[L]} G_n(t,x) q^{L(0)} =
({\rm tr}|_{M(1)}G_n(t,x)q^{L(0)}) ({\rm tr}_{\mathbb{C}[L]} G_n(t,x)q^{L(0)}).
\end{equation}
We fix a basis for $M(1)$:
$$\frac{h(-n_1)^{k_1} \cdots h(-n_l)^{k_l}{\bf 1}}{\sqrt{\prod_{i=1}^l k_i!
n_i^{k_i}}},$$
where $h=\frac{\alpha}{\sqrt{2}}$, $n_1>n_2>...>n_l$
and $k_1,...,k_l \in \mathbb{N}$.
Also we fix the corresponding dual basis
$$\frac{h'(-n_1)^{k_1} \cdots h'(-n_l)^{k_l}{\bf 1}}{\sqrt{\prod_{i=1}^l k_i!
n_i^{k_i}}}.$$ Now we calculate (as in \cite{ts2}) \bea && {\rm
tr}_{\mathbb{C}[h(-m)]}G_n(t,x)q^{L(0)}= \nn && \sum_{k \in
\mathbb{N}} \frac{1}{k! m^k} \langle h'(-m)^{k}, G_n(t,x)
q^{L(0)}h(-m)^{k} \rangle =\nn &&\sum_{m \in \mathbb{N}}
\frac{\displaystyle{\prod_{i <j}}
((t_jx_j/t_ix_i)^{1/2}-(t_jx_j/t_ix_i)^{-1/2})
((x_j/x_i)^{1/2}-(x_j/x_i)^{-1/2})}{\displaystyle{\prod_{i=1}^n}(t_i^{1/2}-t_i^{-1/2})
\displaystyle{\prod_{i <j}}
((t_jx_j/x_i)^{1/2}-(t_jx_j/x_i)^{-1/2})((x_j/t_ix_i)^{1/2}-(x_j/t_ix_i)^{-1/2})
} \nn && \frac{1}{k! m^k} \langle h'(-m)^{k} \nb G_n(t,x) \nb
h(-m)^{k} \rangle ,\nonumber. \eea Then \bea && \sum_{k \in
\mathbb{N}}\frac{1}{k! m^k} \langle h'(-m)^{k} \nb G_n(t,x) \nb
h(-m)^{k} \rangle =\nn && \sum_{k \in \mathbb{N}}\frac{1}{k! m^k}
\langle h'(-m)^{k}, {\rm
exp}\left(\frac{\alpha(-m)}{m}\sum_{r=1}^n
(x_st_s)^{-m}+x_s^{-m}\right) \nn &&{\rm
exp}\left(\frac{\alpha(m)}{-m}\sum_{i=1}^n
(x_rt_r)^{nm}+x_r^{m}\right) h(-m)^{k} \rangle= \nn &&\sum_{k \in
\mathbb{N}} \frac{1}{k! m^k} \sum_{i=0}^k \frac{1}{i!}
(-\frac{1}{m} \sum_r (x_rt_r)^{-m}+x_r^{-m})^i  \langle
h'(-m)^{k}, {\rm exp} \biggl(\frac{\alpha(-m)}{m} \nn &&
\sum_{s=1}^n (x_st_s)^{-m}+x_s^{-m}\biggr) \alpha(m)^i q^{mk}
h(-m)^{k} \rangle =\nn && \sum_{i=0}^\infty \sum_{k=i}^\infty
\frac{k! \langle h'(-m)^k,h(-m)^k \rangle }{m^k k!}
\frac{(-\frac{1}{m}\sum_{r,s} \frac{t_r^mx_r^m} {t_s^m x_s^m}+
\frac{x_r}{x_s}-\frac{t_r^mx_r^m}{x_s}-\frac{x_r}{t_s^mx_s^m})^i}{(k-i)!
i!^2}q^{mk}=\nn &&\sum_{i=0}^\infty \frac{(-\frac{1}{m}\sum_{r,s}
\frac{t_r^mx_r^m} {t_s^m
x_s^m}+\frac{x_r}{x_s}-\frac{t_r^mx_r^m}{x_s}-\frac{x_r}{t_s^mx_s^m})^i\frac{q^{im}}{(1-q^m)^i}}{i!(1-q^m)}=\nn
&&{\rm exp}\left(-\frac{1}{m}\sum_{r,s} (\frac{t_r^mx_r^m} {t_s^m
x_s^m}+
\frac{x_r}{x_s}-\frac{t_r^mx_r^m}{x_s}-\frac{x_r}{t_s^mx_s^m})\frac{q^m}{1-q^m}\right)
\frac{1}{1-q^m}. \eea So far we calculated the trace on the space
$\mathbb{C}[h(-m)]$. Since
$$M(1) \cong \otimes_{m=1}^{\infty} \mathbb{C}[h(-m)],$$
then
$${\rm tr}|_{M(1)} A=\prod_{m=1}^\infty tr|_{\mathbb{C}[h(-m)]}A. $$
Therefore
\bea \label{infiniteprod}
&& \prod_{m=1}^\infty  exp\left(-\frac{1}{m}(\frac{t_r^mx_r^m}
{t_s^m x_s^m}+\frac{t_s^mx_s^m}{t_r^m
x_r^m}+1)\frac{q^m}{1-q^m}\right)= \\
&& \prod_{t=1}^{\infty} \frac{(1-q^t \frac{t_r x_r}{t_s x_s})(1-q^t
\frac{t_s x_s}{t_r x_r})}{(1-q^t)^2}, \nonumber
\eea
for  $1 \leq r \leq n$ and $1 \leq s \leq n$.
Thus
\bea
&&\frac{\prod_{i <j} ((t_jx_j/t_ix_i)^{1/2}-(t_jx_j/t_ix_i)^{-1/2})
((x_j/x_i)^{1/2}-(x_j/x_i)^{-1/2})}{\prod_{i=1}^n(t_i^{1/2}-t_i^{-1/2})
\prod_{i <j} ((t_jx_j/x_i)^{1/2}-(t_jx_j/x_i)^{-1/2})((x_j/t_ix_i)^{1/2}-(x_j/t_ix_i)^{-1/2})
} \nn
&& \prod_{m=1}^\infty exp\left(-\frac{1}{m}\sum_{r,s} (\frac{t_r^mx_r^m}
{t_s^m x_s^m}+
\frac{x_r}{x_s}-\frac{t_r^mx_r^m}{x_s}-\frac{x_r}{t_s^mx_s^m})\frac{q^m}{1-q^m}\right)
\frac{1}{1-q^m}= \nn
&&  \frac{\prod_{i <j} \theta_{11}(t_j x_j/t_ix_i) \theta_{11}(x_j/x_i)}
{\prod_{i<j} \theta_{11}(t_j x_j/x_i) \theta_{11}(x_j/t_ix_i)
\prod_{i=1}^n \theta_{11}(t_i) (q)_{\infty}}.
\eea
If we assume that
$$|t_1x_1|>|x_1|>...>|t_nx_n|>|x_n|>|qt_1x_1|>0,$$
then the infinite
product (\ref{infiniteprod}) converges uniformly
to a holomorphic function inside
(\ref{domain2}). This is a proof of (b).
For (a) we have to calculate (cf. (\ref{prodtr}) )
$${\rm tr}|_{\mathbb{C}[L]} G_n(t,x)q^{L(0)}=\sum_{n \in \mathbb{Z}}
t_1^m \cdots t_n^m q^{\frac{m^2}{2}}.$$
which is convergent for $|q|<1$.
\epfv

From Theorem \ref{mother} it follows that
$$\eta(q) {\rm tr}|_{{\cal F}_0}
X(\psi,t_1x_1)X(\psi^* ,x_1) \cdots X(\psi,t_nx_n)X(\psi^* ,x_n)
q^{L(0)}$$ has a meromorphic continuation (after we specify
branches) to $(\mathbb{C}^\times)^{2n}$ with the set of
singularities determined by zeros of \be \label{singularities}
\prod_{i=1}^n \theta_{11}(t_i) \prod_{i<j} \theta_{11}(t_j
x_j/x_i) \theta_{11}(x_j/t_ix_i). \ee Let $x=(x_1,...,x_n)$.
Denote by $\Lambda_n(x)$ the set of all $n$--tuples
$(t_1,...,t_n)$ such that $(t_1,...,t_n,x_1,...,x_n)$ satisfies
inequality (\ref{domain2}).


\subsection{Another description of the $q$--difference equations}

From Theorem \ref{mother}(b)
it follows:
\begin{corollary} \label{onetheta}
Let $|q|<1$. Then
\begin{equation} \label{11ll}
{\rm tr}|_{\mathcal{F}_0}
 \eta(q) X(a,t x)X(a^{-1},x) q^{L(0)}
\end{equation}
converges in the domain $\frac{1}{|q|}>|t|>1$.
Moreover, (\ref{11ll}) has a meromorphic extension to
a double cover of $\mathbb{C}^\times$ equal to
$$\frac{1}{\theta_{11}(t)}.$$
\end{corollary}
Corollary \ref{onetheta} is essentially Theorem 6.5 in
\cite{BO}.
\begin{remark} \label{disc}
{\em According to \cite{BO}, the following formula holds:
\begin{equation} \label{blochoku}
{\cal F}(y_1,...,y_n)=\sum_{\sigma \in S_n}
\frac{{\rm det} \left(
\frac{\theta_{11}^{(i-j+1)} ( t_{\sigma(1)} \cdots t_{\sigma(n-j)} )}{(j-i+1)!}
\right)_{i,j=1}^n}
{\theta_{11} ( t_{\sigma(1)} )...\theta_{11}(
t_{\sigma(1)}...t_{\sigma(n)} ) },
\end{equation}
where $t_i=e^{y_i}$.
Let us denote by $G(t_1,...,t_n)$ the coefficient of $x^0$ in the
expansion of (\ref{ttq}).
Because of Proposition \ref{main3}
$${\cal F}(y_1,...,y_n)=G(t_1,...,t_n)|_{t_i^{n/2}=e^{n y_i/2}}.$$
If we switch to complex variables
a Laurent expansion of (\ref{thetaproducts}) (inside a certain domain)
gives us an integral representation of $G(t_1,...,t_n)$.
Let $A_i$ be an annulus inside the $x_i$--plane
such that $A_i$ is contained in (\ref{domain2}), and
${\cal C}_i$ is a smooth curve inside $A_i$. Then
the Cauchy's integral formula for several complex variables implies
\begin{equation} \label{blochokuint}
{G}(t_1,...,t_n)=\frac{1}{(2 \pi i)^n}
\oint_{C_1} \cdots \oint_{C_n} \frac{dx_1 \cdots dx_n}{x_1
\ldots x_n} \frac{\displaystyle{\prod_{i <j}} \theta_{11}(t_j x_j/t_ix_i) \theta_{11}(x_j/x_i)}
{\displaystyle{\prod_{i<j}} \theta_{11}(t_j x_j/x_i) \theta_{11}(x_j/t_ix_i)
\displaystyle{\prod_{i=1}^n} \theta_{11}(t_i)},
\end{equation}
This integral is hard to calculate explicitly.}
\end{remark}

Let us consider the case $n=2$ in more details.
Let $x=\frac{x_2}{x_1}$. Assume that $q^n=e^{2 \pi n \tau}$, $\tau
\in \mathbb{H}$ and $|t_2|>|t_1|$. Then the set of singularities (\ref{singularities})
is ordered (with respect to the absolute value) in the
following way:
\begin{equation} \label{order}
\ldots > \left|\frac{q^k}{t_2}\right|>|t_1|>\left|\frac{q^{k+1}}{t_2}\right|>|q
t_1|>\left|\frac{q^{k+2}}{t_2}\right|> \ldots,
\end{equation}
for some $k \in \mathbb{Z}$. Suppose for simplicity that $k=0$, i.e.
\begin{equation} \label{order1}
\ldots >
\left|\frac{1}{t_2}\right|>|t_1|>\left|\frac{q}{t_2}\right|>|q t_1|>\left|\frac{q^{2}}{ t_2}\right|> \ldots.
\end{equation}
Put $r_1=|q t_1|$, $r_2=\left|\frac{q}{t_2}\right|$ and $r_3=|t_1|$.
Let us denote by ${\cal C}_r(z) \in \mathbb{C}$ a circle with center at $z$ and radius $r$.
From Remark \ref{disc}, in the case when
$n=2$, it follows
\begin{corollary}
Let
\begin{equation} \label{domain2pt}
|t_1|>\left|\frac{q}{t_2}\right|.
\end{equation}
Then
\bea \label{2elint}
{G}(t_1,t_2)= \frac{1}{\theta_{11}(t_1)\theta_{11}(t_2)} \oint_{{\cal C}_{a_1}(0)}
\frac{dx}{2 \pi i x} \frac{\theta_{11}(t_2 x/t_1) \theta_{11}(x)}
{\theta_{11}(t_2 x) \theta_{11}(x/t_1)},
\eea
where
$$a_1=\frac{|t_1|+\frac{|q|}{|t_2|}}{2}.$$
\end{corollary}
Let us try to determine the behavior of (\ref{2elint}) under the
elliptic transformation $t_1 \mapsto qt_1$.
Now,
$$G(qt_1,t_2)=\frac{1}{\theta_{11}(q t_1)\theta_{11}(t_2)} \oint_{{\cal C}_{a_2}(0)}
\frac{dx}{2 \pi i x} \frac{\theta_{11}(t_2 x/q t_1) \theta_{11}(x)}
{\theta_{11}(t_2 x) \theta_{11}(x/q t_1)},$$
where
$$a_2=\frac{\frac{|q|}{|t_2|}+|q t_1|}{2}.$$
Since
$$G(qt_1,t_2,x)=-q^{1/2}t_1t_2 G(t_1,t_2,x),$$
$$G(qt_1,t_2)=-\frac{q^{1/2}t_1t_2}{\theta_{11}(t_1)\theta_{11}(t_2)} \oint_{{\cal C}_{a_2}(0)}
\frac{dx}{2 \pi i x} \frac{\theta_{11}(t_2 x/ t_1) \theta_{11}(x)}
{\theta_{11}(t_2 x) \theta_{11}(x/ t_1)},$$
\begin{center}
\begin{figure}
{\par \centering \resizebox*{!}{5cm}{\includegraphics{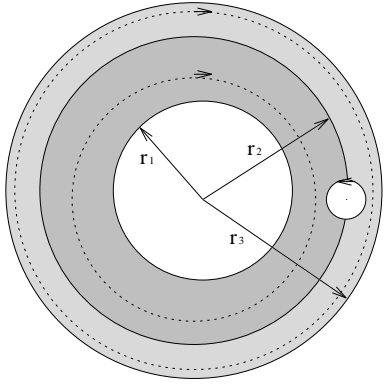}} \par}
\caption{\label{circle} Contours of integration}
\end{figure}
\end{center}


Now from (\ref{order1})
and Cauchy's theorem
we have
\bea \label{homology}
&& G(qt_1,t_2)= -\frac{q^{1/2}t_1t_2}{\theta_{11}(t_1)\theta_{11}(t_2)} \oint_{{\cal C}_{a_2}(0)}
\frac{dx}{2 \pi i x} \frac{\theta_{11}(t_2 x/ t_1) \theta_{11}(x)}
{\theta_{11}(t_2 x) \theta_{11}(x/ t_1)}=\nn
&& -q^{1/2}t_1t_2 \biggl( \frac{1}{\theta_{11}(t_1)\theta_{11}(t_2)}
\oint_{{\cal C}_{a_1}(0)}
\frac{dx}{2 \pi i x} \frac{\theta_{11}(t_2 x/ t_1) \theta_{11}(x)}
{\theta_{11}(t_2 x) \theta_{11}(x/ t_1)}
+ \nn
&&  \frac{1}{\theta_{11}(t_1)\theta_{11}(t_2)} \oint_{{\cal C}_{\epsilon}(\frac{q}{t_2})}
\frac{dx}{2 \pi i x} \frac{\theta_{11}(t_2 x/ t_1) \theta_{11}(x)}
{\theta_{11}(t_2 x) \theta_{11}(x/t_1)} \biggr),\nn
\eea
for some small $\epsilon>0$  (cf. Figure \ref{circle}, where
dotted circles have radius $a_1$ and $a_2$ and the small circle
is centered at $\frac{q}{t_2}$).
Let us evaluate the second integral on the right hand side of (\ref{homology}).

\begin{lemma} \label{defcont}
For small $\epsilon >0$
\begin{itemize}
\item[(a)]
\begin{equation} \label{third}
\frac{1}{\theta_{11}(t_1)\theta_{11}(t_2)} \oint_{{\cal C}_{\epsilon}(qt_1)}
\frac{dx}{2 \pi i x} \frac{\theta_{11}(t_2 x/ t_1) \theta_{11}(x)}
{\theta_{11}(t_2 x) \theta_{11}(x/t_1)}=\frac{1}{\theta_{11}(t_1t_2)},
\end{equation}
and
\item[(b)]
\begin{equation}
\frac{1}{\theta_{11}(t_1)\theta_{11}(t_2)} \oint_{{\cal C}_{\epsilon}(\frac{q}{t_2})}
\frac{dx}{2 \pi i x} \frac{\theta_{11}(t_2 x/ t_1) \theta_{11}(x)}
{\theta_{11}(t_2 x) \theta_{11}(x/t_1)}=\frac{1}{\theta_{11}(t_1t_2)}.
\end{equation}
\end{itemize}

\end{lemma}
{\em Proof:}
The proof of (b) is trivial once we prove (a).
Note first that the two functions in (\ref{third}) are multivalued. Since we
can extract on both sides a term $\sqrt{t_1t_2}$ we may think of
(\ref{third})  as an equality of single--valued functions.
\bea
&& \frac{1}{\theta_{11}(t_1)\theta_{11}(t_2)} \oint_{{\cal C}_3}
\frac{dx}{2 \pi i x} \frac{\theta_{11}(t_2 x/ t_1) \theta_{11}(x)}
{\theta_{11}(t_2 x) \theta_{11}(x/t_1)}= \nn
&& \frac{1}{\theta_{11}(t_1) \theta_{11}(t_2) } \cdot \nn
&& \oint_{{\cal C}_3}
\frac{dx}{2 \pi i x (1-\frac{qt_1}{x}) } \frac{ \theta_{11}(t_2 x/ t_1) \theta_{11}(x)(q)_{\infty}^2}
{\theta_{11}(t_2 x) \prod_{i \geq 1}
(1-\frac{q^{i+1}t_1}{x}) (1-\frac{q^{i}x}{t_1})
( (x/t_1)^{1/2}-(x/t_1)^{-1/2} ) }=\nn
&& \frac{\theta_{11}(qt_2) \theta_{11}(qt_1) (q)_{\infty}^2}
{\theta_{11}(t_1) \theta_{11}(t_2) \theta_{11}(qt_1t_2) \prod_{i \geq
1}(1-q^{i+1})(1-q^i) (q^{1/2}-q^{-1/2}) }=\nn
&& \frac{\theta_{11}(t_2) \theta_{11}(t_1) (q)_{\infty}^2} {\theta_{11}(t_1)\theta_{11}(t_2)
\theta_{11}(t_1t_2)(q)_{\infty}^2}=\frac{1}{\theta_{11}(t_1t_2)}.
\eea
\epfv

Lemma \ref{defcont} and (\ref{homology}) imply that
the $q$--difference equation for $G(t_1,t_2)$ depends on the
one--point function $G(t_1 t_2)$, i.e.
$$G(qt_1,t_2)=-q^{1/2}t_1 t_2 (G(t_1,t_2)-G(t_1t_2)).$$
In the next section we will prove the same formula
by using purely algebraic methods.

\subsection{$q$--difference equations for Bloch-Okounkov $n$--point
functions}

In this part we analyze the difference
equations for ${ G}(t_1,...,t_n)$ (at least when $n=2$)
which gives us elliptic transformation formulas for the
multi valued function ${G}(t_1,...,t_n)$,
i.e. transformation formulas for ${G}(q t_1,...,t_n)$.
If we put
$$H(t_1,...,t_n,x_1,...,x_n)=\frac{\ds{\prod_{i <j}}
\theta_{11}(t_j x_j/t_ix_i) \theta_{11}(x_j/x_i)}
{\ds{\prod_{i<j}} \theta_{11}(t_j x_j/x_i)
\theta_{11}(x_j/t_ix_i) \ds{\prod_{i=1}^n } \theta_{11}(t_i) },$$
then we have
\begin{equation} \label{ellipticbulk}
H(qt_1,...,t_n,x_1,...,x_n)=-q^{1/2}t_1t_2...t_n
H(t_1,...,t_n,x_1,...,x_n).
\end{equation}
where we used the  formula
$$\theta_{11}(qt)=-q^{-1/2}t^{-1}\theta_{11}(t).$$
We are aiming for a similar formula
as in the case of $G(t_1,...,t_n)$.
For simplicity let us consider first
the $1$ and $2$--point functions.

We start with the $1$--point function
\bea \label{brr}
{ G}(t_1)=\eta(q) {\rm tr}|_{{\mathcal{F}}_0} \ X(\psi,t_1x_1)X(\psi^*,x_1) q^{\bar{L}(0)}.
\eea
Then we have
\begin{proposition} \label{oneq}
\begin{equation} \label{forone}
{ G}(q t_1)+q^{1/2}t_1 { G}(t_1)=\delta_{1/2}(q t_1).
\end{equation}
\end{proposition}
{\em Proof:}
First note that
$$[X(\psi,t_1x_1), X(\psi^*,x_1)]=\delta_{1/2}( t_1).$$
Also, we introduce a ``control'' variable $y$ a $\mathbb{C}[y,y^{-1}]$--valued operator
$y^{H(0)}$, where $H(0)$ is the charge operator introduced earlier.
\bea \label{recursion1}
&& \eta(q) {\rm tr}|_{\cal{F}} X(\psi,t_1x_1)X(\psi^*,x_1)y^{H(0)} q^{\bar{L}(0)}= \nn
&& \eta(q) {\rm tr}|_{\cal{F}} \delta_{1/2}( t_1) y^{H(0)} q^{\bar{L}(0)}-
\eta(q) {\rm tr}|_{\cal{F}} X(\psi^*,x_1)X(\psi,t_1x_1)y^{H(0)} q^{\bar{L}(0)}= \nn
&&\eta(q) {\rm tr}|_{\cal{F}} \delta_{1/2}( t_1)q^{\bar{L}(0)}-{\rm
tr}|_{\cal{F}}\eta(q)
X(\psi^*,x_1)y^{H(0)-1}q^{\bar{L}(0)}X(\psi,\frac{t_1x_1}{q})=\nn
&& \eta(q) {\rm tr}|_{\cal{F}}  \delta_{1/2}( t_1)q^{\bar{L}(0)}- \eta(q) {\rm
tr}|_{\cal{F}}
X(\psi,\frac{x_1t_1}{q})X(\psi^*,x_1)y^{H(0)-1}q^{\bar{L}(0)}.
\eea
Now if we extract the coefficient of  $y^{0}$ in (\ref{recursion1})
we obtain
$${ G}(t_1)=\delta_{1/2} (t_1)-q^{-1/2}t_1{G}(t_1/q).$$
After multiplying with $q^{D_{t_1}}$ we obtain the desired formula.
\epfv

Note that  (\ref{forone}) is a formal analogue of the difference equation
for one-point functions obtained in \cite{BO}.
By using the simple identity
$(1-t)\delta_{1/2}(t)=0$.
we obtain
$$(1-qt_1)({G}(q t_1)+q^{1/2}t_1 {G}(t_1))=0.$$
Note that in the equation above we cannot remove the term
$(1-qt_1)$ because this is a formal expression. On the other hand
this equation gives us meromorphic continuation on (a double cover
of) $\mathbb{C}^\times$ (provided that ${G}(t_1)$ is
meromorphic in $\Omega$ and such that $q^\mathbb{Z}
\cdot \Omega=\mathbb{C}^\times$). Now we consider transition from the
formal variables to the complex numbers. From (\ref{forone}) it
follows \bea && {G}(q t_1)+q^{1/2}t_1 {
G}(t_1)=(q^{D_{t_1}}+q^{1/2}t_1){G}(t_1)= \nn &&
=(1+q^{D_{t_1}}q^{-1/2}t_1^{-1})q^{1/2}t_1{G}(t_1) = \nn &&=
\delta_{1/2}(q t_1). \eea Now we formally invert the operator
$(1+q^{D_{t_1}}q^{-1/2}t_1^{-1})$ by using a formal expansion
$$\frac{1}{1+A}=\sum_{n \geq 0}(-1)^n A^n.$$
Then we have
\bea
&& q^{1/2}t_1{G}(t_1)=\frac{1}{1+q^{D_{t_1}}q^{-1/2}t_1^{-1}}\delta_{1/2}(q
t_1)=\nn
&& =\frac{1}{1+q^{-3/2}t_1^{-1}q^{D_{t_1}}}\delta_{1/2}(q
t_1).
\eea
Therefore
\bea
&&{G}(t_1)=q^{-1/2}t_1^{-1} \sum_{m \geq 0} (-1)^m
(q^{-3/2}t_1^{-1}q^{D_{t_1}})^m \delta_{1/2}(q t_1)=\nn
&&=q^{-1/2}t_1^{-1} \sum_{n \in {\mathbb Z}+1/2}\sum_{m \geq 0}
(q^{-3/2}t_1^{-1}q^{D_{t_1}})^m q^n t_1^n=\nn
&& \sum_{k \in {\mathbb Z}}\sum_{m \geq 0}
(-1)^m q^{\frac{m(m+1)}{2}}q^{-k(m+1)}t_1^{-k-1/2},
\eea
where $n-m=-k+1/2$.
Now, since
$$\frac{m(m+1)}{2}-k(m+1)=\frac{(m-k+1/2)^2}{2}-\frac{(k+1/2)^2}{2},$$
 we obtain
\begin{equation} \label{almtheta}
{ G}(t_1)=\sum_{k \in {\mathbb Z}} t_1^{-k-1/2}
q^{-\frac{(k+1/2)^2}{2}}
 \sum_{m \geq 0} (-1)^m q^{\frac{(m-k+1/2)^2}{2}}.
\end{equation}
By using the symmetry on the right hand side of
 (\ref{almtheta}) and after summing
the geometric series we obtain:
\begin{equation} \label{conv}
G(t_1)=t_1^{-1/2}\sum_{m \geq 0} (-1)^m q^{\frac{m(m+1)}{2}}
\left(\frac{q^m t_1^{-1}}{1-t_1^{-1}q^m}+1+
\frac{q^{m+1} t_1}{1-t_1 q^{m+1}}\right).
\end{equation}
If we switch to complex variables then the
Note that (\ref{conv}) is a convergent series when
$1<|t|< \frac{1}{|q|}$, and it can be analytically extended
(by using (\ref{conv})) to all values $t \neq q^n$, $n \in \mathbb{Z}$.
Thus
\begin{equation} \label{sumtheta}
{G}(t_1)=t_1^{-1/2} \sum_{m \geq 0} (-1)^m q^{\frac{m(m+1)}{2}}
\frac{1-q^{2m+1}}{(1-t_1^{-1}q^m)(1-t_1 q^{m+1})}.
\end{equation}
We claim that ${G}(t_1)$ in (\ref{sumtheta}) is
equal to (\ref{brr}) (remember we just find a solution of
(\ref{oneq})). For this apply Lemma \ref{appl2} from
Appendix A.

Now, if we combine formula (\ref{sumtheta}) with Corollary
(\ref{onetheta}) we see that the reciprocal theta function can be
written as an infinite sum (with the same set of poles).
The similar formula follows from the denominator formula for $N=2$
superconformal algebra.

Now we discuss two point functions. Again the aim is to derive
q--difference equation for the two point function by using a formal
variable approach.
\begin{proposition} \label{twoq}
We have
\begin{equation} \label{twopt}
G(qt_1,t_2)+q^{1/2} t_1t_2 G(t_1,t_2)=\delta_{1/2}(q
t_1)G(t_2)+q^{1/2}t_1 t_2 G(t_1t_2).
\end{equation}
Moreover we have a locality equation:
\begin{equation} \label{twoptloc}
(1-q t_1 t_2)(1-q t_1) \left\{ G(qt_1,t_2)+q^{1/2}t_1t_2G(t_1,t_2)-q^{1/2}t_1
t_2G(t_1t_2) \right\}=0.
\end{equation}
\end{proposition}
{\em Proof:}
The proof is very similar as the one in Proposition (\ref{oneq}).
So we omit details.  By using the Jacobi identity
\bea \label{2pointqdiff}
&& {\rm tr}|_{\cal{F}}\eta(q)o(X(\psi,t_1x_1)X(\psi^*,x_1))o(X(\psi,t_2 x_2)X(\psi^*,x_2))
y^{H(0)} q^{\bar{L}(0)}+ \nn
&& {\rm tr}|_{\cal{F}}\eta(q)o(X(\psi,\frac{t_1x_1}{q})
X(\psi^*,x_1))o(X(\psi,t_2 x_2)X(\psi^*,x_2))
y^{H(0)-1} q^{\bar{L}(0)}=\nn
&& {\rm tr}|_{\cal{F}}\eta(q)\delta_{1/2}(t_1)
X(\psi,t_2x_2)X(\psi^*,x_2)y^{H(0)} q^{\bar{L}(0)}-\nn
&& {\rm tr}|_{\cal{F}} X(\psi, t_1t_2 x_1) X(\psi^*,x_1)
y^{H(0)} q^{\bar{L}(0)}+{\rm tr}|_{\cal{F}} \delta_{1/2}(t_1 t_2)
y^{H(0)} q^{\bar{L}(0)}.
\eea
If we extract coefficient of $(x_1 x_2)^0$ in (\ref{2pointqdiff}), by using
Theorem 3.1, we obtain
\bea \label{2qeq}
&& {G}(t_1,t_2)+q^{-1/2}t_1 t_2  {G}\left(\frac{t_1}{q}, t_2
\right) =\nn
&& \delta_{1/2}(t_1) G(t_2)-G (t_1 t_2 )+\delta_{1/2}(t_1 t_2).
\eea
If we act by $q^{D_t}$ we obtain
\bea
&& {G}(q t_1,t_2)+q^{1/2}t_1 t_2  {G}\left(t_1, t_2
\right) =\nn
&& q^{1/2}t_1 t_2 G(t_1 t_2)-\delta_{1/2}(q t_1 t_2)+\delta_{1/2}(q t_1)
G(t_2)+\delta_{1/2}(q t_1 t_2)=\nn
&& q^{1/2}t_1 t_2 G(t_1 t_2)+\delta_{1/2}(q t_1)
G(t_2) .
\eea
\epfv

By using this method we can derive similar formulas for
an arbitrary $n$--point function. A similar theorem
has been proven in  \cite{BO} by using different methods.
\begin{theorem} \label{mainmain}
Formally
\begin{itemize}
\item[(a)]
\bea \label{canada}
&& G(t_1,...,t_n)+q^{-1/2} \left( \prod_{i=1}^n t_i \right) G
\left(\frac{t_1}{q},t_2,...,t_n \right)=\nn
&& - \sum_{s=1}^{n-1} \sum_{1< t_{i_1}< \ldots <t_{i_{s}} \leq n  }
G(t_1t_{i_1} \cdots t_{i_{s}},
\ldots,\hat{t}_{i_1},\ldots, \hat{t}_{i_k},\ldots )+\nn
&&+ \delta_{1/2}(t_{i_1} \cdots t_{i_s}) G(\cdots \hat{t}_{i_1} \cdots \hat{t}_{i_s} \cdots).
+ \delta_{1/2}(t_1)G(t_2,...,t_n).
\eea
\item[(b)]
\bea \label{canada-1}
&& \prod_{0 \leq s \leq n-1, 1<i_1<...i_s \leq n} (1-q t_1 t_{i_1} \ldots t_{i_s}) \biggl(
G(q t_1,...,t_n)+ \\
&& q^{1/2}\prod_{i=1}^n  t_i
\sum_{s=0}^{n-1} \sum_{1<i_1<...<i_s \leq n} (-1)^{s} G(t_1 t_{i_1} \cdots
t_{i_s},\ldots \hat{t}_{i_1}, \ldots ,\hat{t}_{i_s}, \ldots) \biggr)=0 \nonumber.
\eea
\item[(c)]
Let $|q|<1$. Then $G(t_1,...,t_n)$ converges inside an open subset
$$\Omega_{n,0} \cap \{ (t_1,...,t_n): |\prod_{i=1}^n t_i|<\frac{1}{q} \},$$
and it has a meromorphic continuation on a double cover of
$(\mathbb{C}^\times)^n$ with the set of poles contained in
$\{(t_1,...,t_n): q^m t_{i_1} \ldots t_{i_n}=1, 1 \leq i_1<...<i_s
\leq n, m \in \mathbb{Z} \}$.
\end{itemize}
\end{theorem}
{\em Proof:}
The idea is essentially the same as in (\ref{2pointqdiff}).
We prove (a) directly and use the induction for (b).
Actually (b) follows if
\bea \label{canada-2}
&& G(q t_1,...,t_n)+ q^{1/2}\prod_{i=1}^n t_i G(t_1,...,t_n)= \nn
&& q^{1/2}\prod_{i=1}^n t_i
\sum_{s=1}^{n-1} \sum_{1<i_1<...<i_s \leq n} (-1)^{s-1} G(t_1 t_{i_1} \cdots
t_{i_s},\ldots ,\hat{t}_{i_1}, \ldots ,\hat{t}_{i_s} \ldots) \nn
&& + \sum_k \sum_{j_1,...,j_k} \alpha_{j_1,...,j_k} \delta(q t_{j_1} \cdots
t_{j_k})G(\ldots, \hat{t}_{j_1}, \ldots \hat{t}_{j_k}, \ldots ),
\eea
where $\alpha_{j_1,...,j_k}$ are some constant.
By multiplying (\ref{canada-2}) by
$$\prod_{\stackrel{0 \leq s \leq n-1}{1<i_1<...<i_s \leq n}} (1-q t_1 t_{i_1} \cdots t_{i_s})$$
we get (\ref{canada-2}).

For $n=1$ and $n=2$ (\ref{canada-2}) holds.
Suppose that it holds for every $k \leq n-1$.
Then
\bea
&& \sum_{s=1}^k \sum_{1<i_1<...<i_s \leq k} \alpha_s G(qt_1 t_{i_1} \cdots
t_{i_s},\ldots ,\hat{t}_{i_1}, \ldots ,\hat{t}_{i_s}, \ldots)= \nn
&&  \sum_{r=1}^k  (-1)^{r-1} \sum_{1 <j_1<...<j_r \leq k}
 \sum_{m=1}^r {r \choose m} (-1)^m \alpha_m G(t_1 t_{i_1} \cdots
t_{i_r},\ldots, \hat{t}_{i_1}, \ldots \hat{t}_{i_r}, \ldots)+ \nn
&& \ldots .
\eea
Where dots stand for terms that involve $\delta(q t_{i_1} \cdots
t_{i_k})$.
In particular if for every $i$, $\alpha_i=1$  then
$$\sum_{m=1}^r {r \choose m} (-1)^m \alpha_m=-1,$$
and
\bea \label{canada0}
&& \sum_{s=1}^k \sum_{i_1<...<i_s \leq k} G(qt_1 t_{i_1} \cdots
t_{i_s},\ldots, \hat{t}_{i_1}, \ldots, \hat{t}_{i_s}, \ldots)=\nn
&&-\sum_{r=1}^k  (-1)^{r-1} \sum_{1 <j_1<...<j_r \leq k}
G(t_1 t_{i_1} \cdots t_{i_r},\ldots ,\hat{t}_{i_1}, \ldots, \hat{t}_{i_r}
\ldots)+ \ldots
\eea
Let $o$ stands for the operator ${\rm Coff}_{x_1^0 \cdots x_n^0}$.
\bea
&&{\rm tr}|_{\mathcal{F}} \ o(X(\psi, x_1t_1)X(\psi^*,t_1)\cdots X(\psi,
x_nt_n)X(\psi^*,t_n))q^{\bar{L}(0)}y^{H(0)}=\nn
&&{\rm tr}|_{\mathcal{F}} \ o( [X(\psi, x_1t_1),X(\psi^*,x_1)\cdots X(\psi,
x_nt_n)X(\psi^*,x_n)])q^{\bar{L}(0)}y^{H(0)}-\nn
&& {\rm tr}|_{\mathcal{F}} \ o(X(\psi^*,x_1)\cdots X(\psi,
x_nt_n)X(\psi^*,x_n)X(\psi, x_1t_1))q^{\bar{L}(0)}y^{H(0)}.
\eea
Therefore
\bea \label{canada1}
&&{\rm tr}|_{\mathcal{F}} \ o(X(\psi, x_1t_1)X(\psi^*,x_1)\cdots X(\psi,
t_nx_n)X(\psi^*,x_n))q^{\bar{L}(0)}y^{H(0)}+\nn
&& {\rm tr}|_{\mathcal{F}} \ o(X(\psi, \frac{x_1t_1}{q})X(\psi^*,x_1)\cdots X(\psi,
x_nt_n)X(\psi^*,x_n))q^{\bar{L}(0)}y^{H(0)-1}=\nn
&&{\rm tr}|_{\mathcal{F}} \ o(\delta_{1/2}(t_1)X(\psi, t_2x_2)X(\psi^*,x_2)
\ldots X(\psi, t_nx_n)X(\psi^*,x_n))q^{\bar{L}(0)}y^{H(0)}+\nn
&& {\rm tr}|_{\mathcal{F}} \ o(X(\psi^*,x_1)X(\psi, t_2 x_2) \delta_{1/2}\left(
\frac{t_1 x_1}{x_2} \right) \cdots X(\psi,
t_n x_n)X(\psi^*,x_n))q^{\bar{L}(0)}y^{H(0)}+ \ldots \nn
&&  {\rm tr}|_{\mathcal{F}}  o(X(\psi^*,t_1) \cdots X(\psi,
t_n x_n) \delta_{1/2} \left( \frac{t_1 x_1}{x_n}
\right))q^{\bar{L}(0)}y^{H(0)}.
\eea
Now in (\ref{canada1}) we use
a formal delta function substitution property (cf. \cite{FLM})
\be
X(u,t_2 x_2) \delta_{1/2}\left(\frac{t_1 x_1}{x_2} \right)=
X(u,t_1 t_2 x_1) \delta \left(\frac{t_1 x_1}{x_2} \right)
\ee
This fact was used in Part I \cite{M} as well.
Then
\bea \label{canada2}
&&{\rm tr}|_{\mathcal{F}} \ o(X(\psi, t_1 x_1)X(\psi^*,x_1)\cdots X(\psi,
t_n x_n)X(\psi^*,x_n))q^{\bar{L}(0)}y^{H(0)}+\nn
&& {\rm tr}|_{\mathcal{F}} \ o(X(\psi, \frac{t_1 x_1}{q})X(\psi^*,x_1)\cdots X(\psi,
t_n x_n)X(\psi^*,x_n))q^{\bar{L}(0)}y^{H(0)-1}=\nn
&& {\rm tr}|_{\mathcal{F}} \ o(\delta_{1/2}(t_1)X(\psi, t_2 x_2)X(\psi^*,x_2)
\cdots X(\psi, t_n x_n)X(\psi^*,x_n))q^{\bar{L}(0)}y^{H(0)}+\nn
&& {\rm tr}|_{\mathcal{F}} \ o(X(\psi^*,x_1)X(\psi, t_1 t_2 x_1) \cdots X(\psi,
t_n x_n)X(\psi^*,x_n))q^{\bar{L}(0)}y^{H(0)}+\ldots \nn
&& {\rm tr}|_{\mathcal{F}} \ o(X(\psi^*,x_1)X(\psi, t_2 x_2) \cdots X(\psi,
t_1 t_n  x_1))q^{\bar{L}(0)}y^{H(0)}.
\eea
We move the term $X(\psi,t_1  t_i x_1)$, $i=2,...,n$
to the right hand side of (\ref{canada2})
to the left by anticommuting them with operators
$X(\psi^*, x_j)$, $j \leq i$. This anticommuting produces more
terms, etc. If we repeat procedure of anticommuting
and delta function substitution for each term (except the first term) on the right hand side of (\ref{canada2})
we obtain
\bea \label{canada3}
&& {\rm tr}|_{\mathcal{F}} \ o(X(\psi^*,x_1)X(\psi, t_1 t_2 x_1) \cdots X(\psi,
t_n x_n)X(\psi^*,x_n))q^{\bar{L}(0)}y^{H(0)}=\nn
&& {\rm tr}|_{\mathcal{F}} \ \delta_{1/2}(t_1t_2) o(\cdots X(\psi, t_n
x_n)X(\psi^*,x_n))q^{\bar{L}(0)}y^{H(0)}-\nn
&& {\rm tr}|_{\mathcal{F}} \ o(X(\psi, t_1 t_2 x_1)X(\psi^*,x_1)\cdots X(\psi,
t_n x_n)X(\psi^*,x_n))q^{\bar{L}(0)}y^{H(0)},
\eea
for the second term,
\bea \label{canada4}
&& {\rm tr}|_{\mathcal{F}} \ o(X(\psi^*,x_1)X(\psi, t_2 x_2) \cdots X(\psi,
t_1 t_i  x_1) \cdots X(\psi, t_n
x_n)X(\psi^*,x_n))q^{\bar{L}(0)}y^{H(0)}=\nn
&& - \sum_{r=1}^{i-1}  \sum_{1< t_{s_1}< \ldots <t_{s_r} \leq i}
{\rm tr}|_{\mathcal{F}} \  o(X(\psi, t_{s_1} \cdots t_{s_{r}}t_i x_1)X(\psi^*,x_1) \cdots
\nn
&& \cdots X(\psi, t_n x_n)X(\psi^*,x_n))q^{\bar{L}(0)}y^{H(0)}\nn
&&+ {\rm tr}|_{\mathcal{F}} \ \delta_{1/2}(t_1 \cdots \hat{t}_{s_1} \cdots
\hat{t}_{s_{r}} \cdots t_i) \cdot \nn
&& o(X(\psi, t_{s_1} x_{s_1})X(\psi^*, x_{s_1}) \cdots X(\psi,
t_{s_{r}}x_{s_{r}}) \cdots X(\psi^*, x_{s_{r}}))q^{\bar{L}(0)}y^{H(0)}
\eea
for the $i+1$--th term, etc.
Combining (\ref{canada2}), (\ref{canada3}) and (\ref{canada4})
we obtain
\bea \label{canada5}
&& G(t_1,...,t_n)+q^{-1/2} \prod_{i \geq 1}^n t_i  G(\frac{t_1}{q},t_2,...,t_n)=\nn
&& -\sum_{s=1}^{n-1} \sum_{1< t_{i_1}< \ldots <t_{i_{s}}} G(t_1t_{i_1} \cdots t_{i_{s}},
\ldots,\hat{t}_{i_1},\ldots, \hat{t}_{i_k},\ldots )+\nn
&&+\delta_{1/2}(t_{i_1} \cdots t_{i_s}) G(\ldots, \hat{t}_{i_1}, \ldots \hat{t}_{i_s} \ldots).
\eea
This is exactly (\ref{canada}).
Now if we act by $q^{D_{t_1}}$ on (\ref{canada5}) we obtain
\bea \label{canada6}
&&G(qt_1,...,t_n)+q^{1/2}t_1 \ldots t_n G({t_1},t_2,...,t_n)=\nn
&&- \sum_{s \geq 1} \sum_{1< t_{i_1}< \ldots <t_{i_{s}}} G(q t_1t_{i_1} \ldots t_{i_{s}},
\cdots,\hat{t}_{i_1},\cdots, \hat{t}_{i_k},\cdots )+\nn
&&+\delta (q  t_{i_1} \ldots t_{i_s}) G(\ldots \hat{t}_{i_1} \ldots \hat{t}_{i_s} \ldots).
\eea
If we apply now (\ref{canada0}) and the inductive assumption
we obtain
\bea \label{canada7}
&& G(q t_1,...,t_n)+ q^{1/2}t_1 \ldots t_n G({t_1},t_2,...,t_n)= \nn
&& q^{1/2}\prod_{i=1}^n t_i
\sum_{s \geq 1} \sum_{1<i_1<...<i_s \leq n} (-1)^{s-1} G(t_1 t_{i_1} \ldots
t_{i_s},\ldots \hat{t}_{i_1} \ldots \hat{t}_{i_s} \ldots) \nn
&& + \delta (q  t_{i_1} \ldots t_{i_s}) G(\ldots \hat{t}_{i_1} \ldots \hat{t}_{i_s} \ldots)+\ldots
\eea
which is (\ref{canada-2}).
There are essentially two ways of proving (c). First which uses
an explicit $q$--expansion for $G(t_1,...,t_n)$ and the second
that uses Theorem \ref{mother}. The first method is much harder
though. Let us explain the second proof. The crucial observation
is that
\begin{equation}
\Upsilon_n \cap \Omega_{n,0} \subset \Lambda_n(x)
\end{equation}
holds for every $x$ where
$$\Upsilon_n=\{ (t_1,...,t_n): |\prod_{i=1}^n t_i|<\frac{1}{q} \}.$$
This is easy to prove by induction. Since the $2n$--point function
converges absolutely inside $\Omega_{n,n}$ and it has an expansion
in terms of $x_i$'s, all terms converge absolutely, hence the
zeroth term $G(t_1,...,t_n)$ converges absolutely inside
$\Upsilon_n \cap \Omega_{n,0}$. Analytic extension and information
about poles follows from the part (b) by using the induction.

%
%
\epfv

Now let us try to solve (\ref{2qeq}) for $n=2$. For $n \geq 3$ it
is possible to obtain similar explicit formula. Because, the
expansion in $t_i$'s converges in a tiny domain, i.e., $\Upsilon_n
\cap \Omega_{n,0}$, these formulas are hard to express in a nice
way (meaning, that the poles and domain of convergence are
visible). On the contrary, product formulas have very nice
form (cf. \cite{BO}).
\begin{proposition}
\begin{itemize}
\item[(a)]
(\ref{2qeq}) has the unique solution in ${\mathbb{C}}[[q,t_1^{\pm
1},t_2^{\pm 1}]]$ that satisfy ${G}(t_1,t_2)={G}(t_2,t_1)$
and it is given by
\bea \label{nnn}
&& {G}(t_1,t_2)= \nn
&& (t_1t_2)^{-1/2} \sum_{m \geq 0} (-1)^m q^{\frac{m(m+1)}{2}}
\frac{1-q^{2m+1}}{(1-(t_1 t_2)^{-1}q^m)(1-(t_1 t_2) q^{m+1})}+\nn
&& (t_1t_2)^{-1/2}\sum_{(k,l) \in \mathbb{Z}^2,k \neq l}
\frac{t_1^{-l}t_2^{-k}}{1-q^{|l-k|}}\cdot \nn
\lefteqn{ \biggl(
\sum_{m=0}^\infty (-1)^{m+1} q^{\frac{m(m+1+|{\rm min}(k,l)|)}{2}}-
q^{|l-k|}\sum_{m=0}^\infty (-1)^{m+1}q^{\frac{m(m+1+|{\rm max}(k,l)|)}{2}}
\biggr).}\nonumber
\eea
\item[(b)]
Let $|q|<1$. Then ${G}(t_1,t_2)$ absolutely converges inside
$$\Omega_{2,0} \cap \{(t_1,t_2): |t_1 t_2|<\frac{1}{|q|} \}.$$
Moreover, it has a meromorphic
extension to a double cover of $(\mathbb{C}^\times)^2$ with the set of poles contained in
$$(\mathbb{C}^\times)^2 \backslash \{(t_1,t_2):
t_1 t_2 q^m=1, t_1 q^m=1, t_2 q^m=1, m \in \mathbb{Z} \}.$$
\end{itemize}
\end{proposition}
{\em Proof:} The fact about the uniqueness follows from Lemma
\ref{appl4} from the Appendix A. Proof of part (a) is quite
lengthy, but straightforward and essentially uses the same
procedure as in the case of $1$--point function. Here is the main
idea. Formula (\ref{2qeq}) can be written as
\begin{equation} \label{asa}
q^{1/2}t_1 t_2
G(t_1,t_2)=\frac{1}{1+q^{-3/2}t_1^{-1}t_2^{-1}q^{D_{t_1}}}
\left(\delta_{1/2}(qt_1)G(t_2)+q^{1/2}t_1 t_2 G(t_1 t_2) \right).
\end{equation}
The fraction in the above formula can be expanded as in the $n=1$
case. The rest is manipulations with $q$--series.

We can prove part (b) in two
different ways. The fact about the convergence
of $G(t_1,t_2)$ follows either from the previous theorem (which holds for every
$n$) or by using explicit expression
(\ref{nnn}).
Hence by using this fact we can meromorphically extend
(\ref{2qeq}) to remaining values of $t$ by using
$q$--difference for $G(qt_1,t_2)$
and $G(q^{-1}t_1,t_2)$
and then by using $G(t_1,qt_2)$
and $G(t_1,q^{-1}t_2)$ (which is obtained from (\ref{2qeq}) because
$G(t_1,t_2)=G(t_2,t_1)$).
From the $q$--difference equation
it follows that we can not extend the function to those values for
which either $G(t_1 t_2)$ has poles or $G(t_1)$ and $G(t_2)$ hence
the statement follows.
\epfv

It is worth mentioning that (as in the $n=1$ case) we encounter
certain $q$--series that have singular parts consisting of
single terms in the $q$--expansion. Once we combined several of
these series all singular terms cancel. This should be compared
with some results from Part I \cite{M}. Let us illustrate this
using a particular example. An {\em incomplete} $\theta$--function
is an expression of the form
$$A_k(q)=\sum_{m=0}^{\infty} (-1)^m q^{\frac{(m+1-2k)(m+2)}{2}},$$
for some $k \in \mathbb{Z}$.
It is clear that $A_k(q) \in \mathbb{Z}[[q]]$ for $k \leq 0$.
Also it is not hard to show that
$$A_k(q)=\frac{1}{q^k}+\sum_{m=0}^\infty (-1)^{m+1} q^{\frac{m(m+2k+1)}{2}}
\in \frac{1}{q^k}+\mathbb{Z}[[q]],$$
for $k>0$.
Then
$$A_{{\rm min}(k,l)}(q)-q^{|l-k|}A_{{\rm max}(k,l)}(q) \in \mathbb{Z}[[q]],$$
for every choice of constants $l,k\in \mathbb{Z}$, $k \neq l$.
This fact is used in proving (\ref{nnn}).
Also notice that coefficients of $t_1^{-k-1/2}t_2^{-l-1/2}$
in (\ref{nnn}) are given by
$$A_{k,l}=\sum_{m,n \geq 0}(-1)^{m+n}
q^{\frac{(n+m+1)(n+m+2)}{2}-k(m+1)-l(n+1)} \in \mathbb{Z}[[q]],$$
for every $k,l \in \mathbb{Z}$, $k \neq l$.
Some related $q$--series were studied in \cite{KWak}.

\subsection{General case}

So far we only discussed charge $0$ subspace, i.e.
$\mathcal{F}_0 \subset \mathcal{F}$.
However $\mathcal{F}_m$ is also a $\hat{\mathcal{D}}$--module.
In vertex operator algebra language, $\mathcal{F}_m$ is
an irreducible $\mathcal{F}_0$--module ($M(1) \cong \mathcal{F}_0 \cong L_1$ cf. \cite{AFOQ}).
Let us denote by $H_m(t_1,...,t_n,x_1,...,x_n)$ the restriction of (\ref{zadnji})
to the space $\mathcal{F}_m$ and with $G_m(t_1,...,t_n)$ the corresponding
Bloch--Okounkov $n$--point function.
For every $m \in \mathbb{Z}$ we have
$$H_m(qt_1,...,t_n,x_1,...,x_n)=-q^{n+1/2}t_1t_2 \cdots t_n H_m
(t_1,...,t_n,x_1,...,x_n).$$
Also,
$$G_m(t_1,...,t_n)= q^{n^2/2}t_1^n \cdots t_2^n G(t_1,...,t_n).$$
Hence the ``total'' Bloch--Okounkov $n$--point function obtained
by taking the trace over $\mathcal{F}$, and it is given by
$$\sum_{m \in \mathbb{Z}} q^{m^2/2} \left( \prod_{i=1}^n t_i \right)^m
G(t_1,...,t_n)= \theta(\prod_{i=1}^m t_i) G(t_1,...,t_n).$$

\renewcommand{\theequation}{\thesection.\arabic{equation}}
\setcounter{equation}{0}

\section{$q$--traces for $\mathcal{D}^-$; Majorana fermion case}

\subsection{Iterated $2n$--point functions}

Let us recall some notation from first part \cite{M}. Denote by
$F$ the Fock fermionic space associated to a free (or Majorana)
fermion. Let $\varphi=\varphi(-1/2){\bf 1} \in F$, such that
$$X(\varphi,x)=\sum_{n \in \mathbb{Z}} \varphi_n x^{-n-1/2}$$
and
$$[\varphi_m,\varphi_n]=\delta_{m+n,0}.$$
The aim is to study $q$--graded traces of the form:
\begin{equation} \label{iteratedbulk}
{\rm tr}|_{\mathcal{F}} X(Y[\varphi,y_1]\varphi,x_1) \cdots X(Y[\varphi,y_n]\varphi,x_n) q^{\bar{L}(0)}.
\end{equation}
If we modify Corollary \ref{main2} and Remark \ref{multival}, such
that it applies in the vertex operator superalgebra case it
follows that all the information about (\ref{iteratedbulk}) can be
obtained by studying the following $2n$--point function:
\begin{equation} \label{singlebulk}
{\rm tr}|_F X(\varphi,t_1x_1)X(\varphi,x_1) \cdots
X(\varphi,t_nx_n)X(\varphi,x_n) q^{\bar{L}(0)}.
\end{equation}
Again, once we switch to the complex variables
(\ref{singlebulk}) is multi--valued with respect to
$t_i$ variables.

It is well--known that one can obtain a free fermion
from a pair of charged fermions.
Let $\tilde{\varphi}(n)=\frac{\psi(n)+\psi^*(n)}{\sqrt{2}}$, for $n \in
\mathbb{Z}+\frac{1}{2}$.
Then
$$[\tilde{\varphi}(m),\tilde{\varphi}(n)]=\delta_{m+n,0},$$
hence computing (\ref{iteratedbulk}), up to the character, reduces
to the case of correlation functions
with charged fermions. Hence, it is not hard to see that
the correlation function (\ref{singlebulk}) converges inside
$$|t_1x_1|>|x_1|> \ldots >|x_n|>|qt_1x_1|>0,$$
to a (multivalued) analytic function. Moreover, (\ref{singlebulk})
has a meromorphic continuation to a double cover of
$(\mathbb{C}^\times)^{2n}$.
Another way of proving convergence, is by explicit calculations.


We present the calculations for $n=2$.
\bea
&& {\rm tr}|_F X(\varphi,t_1x_1)X(\varphi,x_1)
X(\varphi,t_2x_2)X(\varphi,x_2) q^{\bar{L}(0)}+ \nn
&& {\rm tr}|_F X(\varphi,t_1x_1)X(\varphi,x_1) \cdots
X(\varphi,t_nx_n)X(\varphi,x_n) q^{\bar{L}(0)}
=\nn
&& {\rm tr}|_F [X(\varphi,t_1x_1), X(\varphi,x_1) \cdots
X(\varphi,t_nx_n)X(\varphi,x_n)] q^{\bar{L}(0)}=\nn
\lefteqn{ {\rm tr}|_F \delta_{1/2} \left(\frac{1}{t_1} \right)
X(\varphi,t_2x_2)X(\varphi,x_2)q^{\bar{L}(0)}-
{\rm tr}|_F \delta_{1/2} \left(\frac{t_2x_2}{t_1x_1}\right)
X(\varphi,x_1)X(\varphi,x_2)q^{\bar{L}(0)}} \nn
&&+ {\rm tr}|_F \delta_{1/2} \left(\frac{x_2}{t_1x_1}\right)
X(\varphi,t_1)X(\varphi, t_2x_2) q^{\bar{L}(0)}.
\eea
Hence we get
\bea
&& {\rm tr}|_F X(\varphi,t_1x_1)X(\varphi,x_1) \cdots
X(\varphi,t_nx_n)X(\varphi,x_n) q^{\bar{L}(0)}=\nn
&& G(t_1)G(t_2) {\rm tr}|_F q^{\bar{L}(0)}
-G(t_1x_1/t_2x_2) G(x_1/x_2) {\rm tr}|_F q^{\bar{L}(0)} \nn
&& + G(t_1x_1/x_2)G(t_1/t_2x_2) {\rm tr}_F q^{\bar{L}(0)},
\eea
where $G(t)$ is defined as in (\ref{almtheta}) and
$${\rm tr}|_F q^{\bar{L}(0)}=\frac{\eta(q)^2}{\eta(q^2) \eta(q^{1/2})}.$$

As before we consider
and the corresponding Bloch--Okounkov $n$--point function
\begin{eqnarray} \label{singlebo}
&& D(t_1,...,t_n):= \\
&& = \frac{\eta(q^2) \eta(q^{1/2})}{\eta(q)^2}
{\rm tr}|_F o(X(\varphi,t_1x_1)X(\varphi,x_1)) \cdots
o(X(\varphi,t_nx_n)X(\varphi,x_n)) q^{\bar{L}(0)}. \nonumber
\end{eqnarray}
\subsection{$q$--difference equations for $n$--point functions}

Here we obtain $q$--difference equations for (\ref{singlebo}).
\begin{theorem}
\begin{itemize}
\item[(a)]
\bea
&& D(t_1,...,t_n)=-D(\frac{t_1}{q},t_2,...,t_n)+
\delta_{1/2}(t_1)D(t_2,...,t_n) \\
&& \sum_{s=1}^{n-1} \sum_{1<i_1 <...<i_s
\leq n} \sum_{\epsilon_{i_1},...,\epsilon_{i_s} \in \{-1,1 \}}
D(t_1 t_{i_1}^{\epsilon_{i_1}}
\cdots t_{i_s}^{\epsilon_{i_s}},...,\hat{t}_{i_1},...
\hat{t}_{i_s},...) \nn
&&+\delta_{1/2}(t_{i_1} \cdots t_{i_s})
D(...,\hat{t}_{i_1},...,\hat{t}_{i_s},...). \nonumber
\eea
\item[(b)]
\bea
&& \prod_{s,i_1,...,i_s; \epsilon_{i_1},...,\epsilon_{i_s} \in \{-1,1 \}}(1-qt_{i_1}^{{\epsilon}_{i_1}} \ldots
t^{{\epsilon}_{i_1}}_{i_s}) \biggl(
D(qt_1,...,t_n)+ \\
&& \sum_{s=0}^{n-1} \sum_{1<i_1,...,i_s \leq n}
\sum_{\epsilon_{i_1},...,\epsilon_{i_s} \in \{ -1,1 \} } D(t_1
t_{i_1}^{\epsilon_{i_1}} \cdots t_{i_s}^{\epsilon_{i_s}},
...,\hat{t}_{i_1},...,\hat{t}_{i_s},...) \biggr)=0.\nonumber \eea
\item[(c)] $D(t_1,...,t_n)$ converges inside an open subset of
$\Omega_{n,0}$, and it has a meromorphic
continuation on a double cover of $(\mathbb{C}^\times)^n$ with the set of poles
$$\{ t^{\epsilon_1}_{i_1} \ldots
t^{\epsilon_s}_{i_s}q^m=1, 1 \leq i_1< \cdots <i_n \leq n \in
\mathbb{Z}_{>0}, m \in \mathbb{Z}, \epsilon_j \in \{-1,1 \} \}.$$
\end{itemize}
\end{theorem}
{\em Proof:} Part (a) is essentially the same as the proof of
Theorem \ref{mainmain}. The only difference are factors
$\epsilon_i \in \{1,-1 \}$. Part (b) follows directly from (a) by
acting with $q^{D_{t_1}}$. Part (c) follows from (a) and a similar
argument as in Theorem \ref{mainmain}.

\epfv

\renewcommand{\theequation}{\thesection.\arabic{equation}}
\setcounter{equation}{0}

\section{$q$--traces for $\hat{\mathcal{D}}^+$; bosonic case}

\subsection{$1$-- and  $2$--point functions}

For $k \geq 2$, let us define a normalized $k$--th Eisenstein series
$$\tilde{G}_{2k}(q)=\frac{-B_{2k}}{(2k)!}+\frac{2}{(2k-1)!}\sum_{n=1}^\infty
\sigma_{2k-1}(n)q^n \in \mathbb{C}[[q]],$$
where $\sigma_k$ is the sum of the $k$--powers of the divisors of $n$ and
$B_{2k}$ a Bernoulli number.
If $q=e^{2 \pi i \tau}$ we will write $\tilde{G}_{2k}(\tau)$ instead of
$\tilde{G}_{2k}(q)$.
Also we define the meromorphic functions
$$\zeta(z,\tau)=\frac{1}{z}+\sum_{\omega \in
\mathbb{Z}+\mathbb{Z}\tau} \frac{1}{z-\omega}+\frac{1}{\omega}+\frac{z}{\omega^2},$$
$$\wp_2(z,\tau)=\frac{1}{z^2}+\sum_{\omega \in
\mathbb{Z}+\mathbb{Z}\tau}
\frac{1}{(z-\omega)^2}-\frac{1}{\omega^2}$$
and
$$\wp_{k+1}(z,\tau)=-\frac{1}{k} \frac{d}{dz} \wp_{k}(z,\tau),$$
for $k \geq 2$.
Hence
$$(-1)^k (k-1)! \wp_k(z,\tau)=\left(\frac{\partial }{\partial z} \right)^{k-2} \wp_2(z,\tau).$$

$\zeta(z,\tau)$ is the Weierstrass zeta--function,
with the property:
$$\zeta(z+\tau,\tau)=\zeta(z,\tau)+2 \zeta\left(\frac{\tau}{2},\tau
\right).$$
For $k \geq 2$, $\wp_k$ is  elliptic functions ;
$\wp_2(z,\tau)$ is usually called the Weierstrass $\wp$--function.
In the annulus $$0<|z|<\ds{{\rm min}_{m,n}}|m+n \tau|,$$ $\wp_k$ has
a Laurent expansion
$$\wp_{k}(z,\tau)=\frac{1}{z^k}+(-1)^k \sum_{n=1}^\infty {2n+1 \choose
k-1} G_{2n+2}(\tau) z^{2n+2-k},$$ where $G_{2n+2}(\tau)$ are
genuine Eisenstein series (not normalized). In parallel with
$\tilde{G}_{2k}(q)$ it is convenient to introduce \be
\label{wenorm} \tilde{\wp_{k}}(z,q)=\frac{1}{z^k}+(-1)^k
\sum_{n=1}^\infty {2n+1 \choose k-1} \tilde{G}_{2n+2}(q)
z^{2n+2-k}. \ee Then, for $k \geq 2$, we have \be \label{relation}
(2 \pi i)^k \tilde{\wp_k}(2 \pi i z, \tau)={\wp_k}(z, \tau). \ee
The Weierstrass $\zeta$--function and $\wp_k$ functions are
closely related to certain Jacobi modular forms. Let us fix
(formal) $P_{k+1}$--series for $k \geq 1$ \bea \label{Pfunction}
&& P_{k+1}(t,q)=\frac{1}{k!} \left(\sum_{n=1}^\infty \frac{n^k
t^n}{1-q^n}+\frac{(-1)^{k+1}n^k t^{-n}q^n}{1-q^n} \right). \eea Then
(\ref{Pfunction}) converges inside $1>|t|>|q|$ if $t$ and $q$ are
complex variables. If we let $t=e^{2 \pi i y}$, $q=e^{2 \pi i
\tau}$ then in the limit (cf. \cite{La}): \bea \label{weier} &&
P_{2}(e^{2 \pi i y},q)=\frac{1}{(2 \pi i)^2} \left( \wp_2
(y,\tau)+G_2(\tau) \right), \nn && P_l(e^{2 \pi i
y},q)=\frac{(-1)^l (l-1)!}{(2 \pi i)^l}\wp_{l}(y,\tau) \eea for $l
\geq 3$. Note that we cannot derive relation (\ref{weier}) purely
by using formal variables.

The aim is to find the transformation properties of
\be \label{hopelast}
\frac{1}{{\rm ch}_M(q)} {\rm tr}|_M
X(u_1,t_1x_1)X(v_1,x_1)X(u_2,t_2x_2)X(v_2, x_2)q^{L(0)}
\ee
under the elliptic transformation:
$$t_1 \mapsto qt_1 \ {\rm and} \   t_2 \mapsto qt_2.$$
(because of the symmetry we consider only $t_1 \mapsto qt_1$).

We will assume that (\ref{condition}) holds, i.e., that for every $i \neq j$
and $n \geq 0$:
$$u_i(n)v_j=c_{u_i,v_j} \delta_{ {\rm wt}(u_i)+{\rm wt}(v_j)-1,n} {\bf 1}$$
Note that this condition implies
$$u_i[n]v_j=c_{u_i,v_j,n}{\bf 1},$$
for every $n \geq 0$, where $c_{u,v,n} \in \mathbb{C}$.

\begin{remark}
{\em
\be \label{rarely1}
u(n)v=c_1 \delta_{ {\rm wt}(u)+{\rm wt}(v)-1,n} {\bf 1}
\ee
and
\be \label{rarely2}
u[n]v=c_2 \delta_{ {\rm wt}(u)+{\rm wt}(v)-1,n} {\bf 1}
\ee
rarely hold simultaneously because of the different grading. More
precisely, if ${\rm wt}(u)$ and ${\rm
wt}(v) \geq 1$
then both (\ref{rarely1}) and (\ref{rarely2}) hold
if and only if ${\rm wt}(u)={\rm wt}(v)=1$.
This can be shown by writing explicitly
$v[n]$ in terms of $v(i)$'s (and vice--versa).}
\end{remark}

Now suppose that ${\rm wt}(u), {\rm
wt}(v) \geq 1$ and that for every $i$
\begin{equation} \label{zero}
o(u_i)|_M =0,
\end{equation}
The iterated $2$--point function and its zeroth term is easy to compute.
From the Jacobi identity for $X$--operators (cf. \cite{Le2} and \cite{M})
$$[X(u,x_1), X(v,x_2)]={\rm Res}_y \delta \left(\frac{e^y x_2}{x_1}
\right) X(Y[u,y]v,x_2)$$
it follows
\bea \label{gen1}
&& {\rm tr}|_M X(u,t_1x_1)X(v,x_1)q^{L(0)}=\nn
&& {\rm tr}|_M X(u,\frac{t_1x_1}{q})X(v,x_1)q^{L(0)}
+\sum_{k \geq 0} \frac{c_{u,v,k}}{k!} {\rm tr}|_M D_{t_1}^k \delta(t_1) q^{L(0)},
\eea
where $D_{t_1}=t_1\frac{d}{dt_1}$.
Also we may assume that the summation in (\ref{gen1}) goes from $k=1$
\footnote{Note that for superalgebras this will not be the case}.
Or equivalently
\bea
&& {\rm tr}|_M X(u,qt_1x_1)X(v,x_1)q^{L(0)}-
{\rm tr}|_M X(u,t_1x_1)X(v,x_1)q^{L(0)}=\nn
&&\sum_{k=1} \frac{c_{u,v,k}}{k!}{\rm tr}|_M (D^k \delta)(qt_1)
q^{L(0)} \nonumber,
\eea
 Thus
if we let $F(t_1)=\frac{1}{{\rm ch}_M(q)} {\rm tr}|_M X(u_1,t_1x_1)X(v_1,x_1)q^{L(0)}$,
then
\begin{equation} \label{firstrec}
F(q t_1)-F(t_1)=\sum_{ k \geq 1} \frac{c_{u,v,k}}{k!} D_{t_1}^k \delta (t_1)
\end{equation}
If we multiply the equation above with $(1-t_1)^{k_1+1}$
we obtain $(1-t_1)^k(F(q t_1)-F(t_1))=0$.
By using  Lemma \ref{appl1} from Appendix A,
it follows that the general solution of (\ref{firstrec}) in the space
$$q^{h} \mathbb{C}[[t_1,t_1^{-1},q]],$$
is of the form
$$F_{part}(t_1)+f(q),$$
where $F_{part}(t_1)$ is some particular solution and $f(q) \in
\mathbb{C}[[q]]$.
But $F(t_1)$ does not involve terms that contain only powers of $q$
hence it follows that
$$F(t_1)=\sum_{k \geq 1} \frac{c_{u,v,k}}{k!} (-1)^k P_{k+1}\left(\frac{1}{t_1}\right).$$
Note that formally we have more than one solution
of (\ref{firstrec}) inside  $q^{h}
\mathbb{C}[[t_1,t_1^{-1},q,q^{-1}]]$ (cf. Appendix A).

The $2$--point functions are more interesting.
\bea \label{gen2}
&& {\rm tr}|_M X(u_1,t_1x_1)X(v_1,x_1)X(u_2,t_2x_2)X(v_2,x_2)q^{L(0)}=\nn
&& {\rm tr}|_M X(u_1,\frac{t_1x_1}{q})X(v_1,x_1)X(u_2,t_2x_2)X(v_2,x_2)q^{L(0)}
+\nn
&&{\rm tr}|_M [X(u_1,t_1x_1),X(v_1,x_1)X(u_2,t_2x_2)X(v_2,x_2)]q^{L(0)}=\nn
&& {\rm tr}|_M X(u_1,\frac{t_1x_1}{q})X(v_1,x_1)X(u_2,t_2x_2)X(v_2,x_2)q^{L(0)}+\nn
&&{\rm tr}|_M \sum_{k_1 \geq 1} \frac{c_{u_1,v_1,k_1}}{k_1!} D_{t_1}^{k_1}\delta \left(t_1 \right)
X(u_2,t_2x_2)X(v_2,x_2) q^{L(0)}+\nn
&&\sum_{k_2 \geq 1} \frac{c_{u_1,u_2,k_2}}{k_2!}  {\rm tr}|_M  D_{x_1}^{k_2}\delta \left(\frac{t_1x_1}{t_2 x_2}\right)
X(v_1,x_1)X(v_2,x_2)q^{L(0)}+\nn
&&\sum_{k_3 \geq 1} \frac{c_{u_1,v_2,k_3}}{k_3!}  {\rm tr}|_M  D_{x_1}^{k_3}\delta \left(\frac{t_1x_1}{x_2}\right)
X(v_1,x_1)X(u_2,x_2t_2)q^{L(0)}=\nn
&& {\rm tr}|_M X(u_1,\frac{t_1x_1}{q})X(v_1,x_1)X(u_2,t_2x_2)X(v_2,x_2)q^{L(0)}+\nn
&&\sum_{k_1 \geq 1} \frac{c_{u_1,v_1,k_1}}{k_1!}
{\rm tr}|_M D_{t_1}^{k_1}\delta \left(t_1 \right) X(u_2,t_2x_2)X(v_2,x_2) q^{L(0)}+\\
&&\sum_{k_2 \geq 1} \frac{c_{u_1,u_2,k_2}}{k_2 !}  {\rm tr}|_M  \sum_{i=0}^{k_2} (-1)^i {k_2 \choose
i}D_{t_1}^{k_2-i} \delta \left(\frac{t_1x_1}{t_2 x_2}\right)
D_{t_2}^{i} X(v_1,\frac{t_2x_2}{t_1})X(v_2,x_2)q^{L(0)}+\nn
&&\sum_{k_3 \geq 1} c_{u_1,v_2, k_3} {\rm tr}|_M  \sum_{i=0}^{k_3} (-1)^i {k_3 \choose i}  D_{t_1}^{k_3-i} \delta
\left(\frac{t_1x_1}{x_2}\right)X(v_1,x_1) D_{t_2}^{i}X
(u_2,t_1t_2x_1)q^{L(0)}, \nonumber
\eea
where $k_i \in \mathbb{N}$, $i=1,2,3$.
Since,
$${\rm tr}|_M X(v_1,x_1)D_{t_2}^{k_3-i} X(u_2,t_1t_2x_1)q^{L(0)}=
{\rm tr}|_M D_{t_2}^{k_3-i}X(u_2,\frac{t_1t_2x_1}{q})X(v_1,x_1) q^{L(0)}$$
we can extract zeroth term in (\ref{gen2}).
For simplicity suppose that $u_i=v_j=u$, $i=1,2$ and let
$k_1=k_2=k_3=k$.
After taking the coefficient of $x^0$ in (\ref{gen2}) we obtain:
\bea \label{gen3}
&& F(t_1,t_2)-F \left(\frac{t_1}{q},t_2 \right)=\sum_{k \geq 1} \frac{c_{u,u,k}}{k!} \biggl( (D^{k}\delta)\left(t_1 \right) F(t_2)+\nn
&& (-1)^{k}D_{t_2}^{k} F\left(\frac{t_2}{t_1}\right)+ (-1)^{k}
D_{t_2}^{k} F \left(\frac{t_1 t_2}{q} \right) \biggr).
\eea
If we act by $q^{D_{t_1}}$ on (\ref{gen3}) we obtain
\bea \label{gen4}
&& F(qt_1,t_2)-F(t_1,t_2)=\sum_{k \geq 1} \frac{c_{u,u,k}}{k!} \biggl( (D^{k}\delta)\left(q t_1 \right) F(t_2)+\nn
&& (-1)^{k} D_{t_2}^{k} F\left(\frac{t_2}{q t_1}\right)+ (-1)^{k} D_{t_2}^{k} F(t_1 t_2)).
\eea
By multiplying (\ref{gen4}) with $(1-t_1)^{K}$, where $K$ is big enough, we obtain {\em
locality} formula:
\bea \label{gen5}
&& (1-t_1)^{K}(F(qt_1,t_2)-F(t_1,t_2))=\nn
&&(1-t_1)^{K}  \sum_{ k \geq 1} \frac{c_{u,u,k}}{k!} \biggl( (-1)^k D_{t_2}^{k} F
\left(\frac{t_2}{t_1q}\right)+(-1)^k D_{t_2}^{k} F(t_1 t_2) \biggr).
\eea
The following proposition gives us an explicit formula for $F(t_1,t_2)$.
\begin{proposition} \label{mother2}
Suppose that
\bea \label{realeq}
&& H_k(qt_1,t_2)-H_k (t_1,t_2)= (D^{k}\delta)\left(q t_1 \right) F(t_2)+\nn
&& D_{t_2}^{k} F \left(\frac{t_2}{q t_1}\right)+D_{t_2}^{k} F(t_1 t_2),
\eea
where $F(t)=P_{k+1}\left(\frac{1}{t_1}\right)$.
Then the general solution of (\ref{realeq}) that satisfies
$H_k(t_1,t_2)=H_k(t_2,t_1)$ is of the form
\begin{itemize}
\item[(a)]
\bea \label{gen6}
\lefteqn{ H_k(t_1,t_2)=f_k(q)+P_{k+1}\left(\frac{1}{t_1},q\right)
P_{k+1}\left(\frac{1}{t_2},q\right)- } \\
&& \left(\sum_{n \geq 1}^{\infty} \frac{n^{2k}q^n
\left(\left(\frac{t_1}{t_2}\right)^n+\left(\frac{t_2}{t_1}\right)^n \right)}{(1-q^n)^2}+\sum_{n \geq 1}^{\infty} \frac{n^{2k}q^n
((t_1t_2)^n+(t_1t_2)^{-n})}{(1-q^n)^2}\right), \nonumber
\eea
where $f_k(q) \in \mathbb{C}[[q]]$. In particular if $f(q)=0$ and $|q|<1$ then
$H_k(t_1,t_2)$
converges  uniformly on compact subsets inside
a domain:
$$|t_1t_2q|<1, \left|\frac{q}{t_1t_2}\right|<1, \left|\frac{t_1}{t_2}q\right|<1,
\left|\frac{t_2}{t_1}q\right|<1, |t_1|>1, |t_2|>1.$$
In particular, it has a meromorphic extension to $(\mathbb{C}^\times)^2$.
\item[(b)]
If  $t_i=e^{2 \pi i y_i}$, $i=1,2$ and $q=e^{ 2 \pi i \tau}$, $\tau \in \mathbb{H}$, then
\begin{equation}
H_1(\tau+y_1,y_2)=H_1(y_1,y_2)+\frac{1}{(2 \pi
i)^2}\left( \wp_{2}(y_1-y_2,\tau)+\wp_{2}(y_1+y_2,\tau)+G_2(\tau) \right)
\end{equation}
and
\bea
&& H_k(\tau+y_1,y_2)=H_k(y_1,y_2)+\frac{1}{(2 \pi
i)^{k+1}}\biggl( \left(\frac{\partial}{\partial y_1} \right)^{k-1} \wp_{2}(y_1-y_2,\tau)+ \nn
&& (-1)^{k+1} \left(\frac{\partial}{\partial y_1} \right)^{k-1} \wp_{2}(y_1+y_2,\tau) \biggr),
\eea
for $k \geq 2$, where $$y_1+y_2 \notin \mathbb{Z} \oplus \mathbb{Z} \tau,
\ \ y_1-y_2 \notin
\mathbb{Z} \oplus \mathbb{Z}\tau.$$
\end{itemize}
\end{proposition}
{\em Proof:}
Proof of (a) requires solving  (\ref{gen3}), or (\ref{gen4}), inside
$q^h \mathbb{C}[[t_1^{\pm 1},t_2^{\pm 1},q]]$.
The Lemma \ref{appl3} gives us description of all solutions. Hence we
have to find a particular one.
This can be done in the following way (cf. \cite{Fe}, \cite{Zh1}). We
expand the left hand side
of (\ref{gen3}) in powers of $t_1$. Then we compare coefficients
by powers $t_1^n$ on both sides, so we can determine
coefficient of $t_1^n$ in $F(t_1,t_2)$ (since it is multiplied by
$(1-q^{-n})$). The details are not very illuminated so we omit them here.
This gives us a way to find a particular solution of the $q$--difference equation (\ref{gen3}).

It is straightforward to see the uniform (and absolute)
convergence of (\ref{gen6}) when $f(q)=0$.
%
If $t_i=e^{2 \pi i y_i}$, (\ref{gen6}) is convergent.
The locality formula (\ref{gen5}) gives us a meromorphic
continuation of (\ref{gen6}) on 
$(\mathbb{C}^\times)^2$ such that for every $n \in \mathbb{Z}$
$$t_1t_2 \neq q^n, t_1/t_2 \neq q^n.$$
\epfv

Now if we apply the Proposition \ref{mother2} we see that $F(t_1,t_2)$
that satisfies (\ref{gen3}) can be expressed as a linear combination
in terms of $H_k(t_1,t_2)$, with $f_k(q)=0$.

\begin{remark}
{\em Note that condition (\ref{condition}) is
crucial for our considerations. Otherwise one cannot
obtain a recursion procedure, i.e. a way to express
$4$--point function by using $2$--point functions.
Zhu's recursion (cf. \cite{Zh1}) would give us, in
general, a way to express general $4$--point function in
terms of $3$--point functions.}
\end{remark}

Explicit calculation of the zeroth term of (\ref{gen2})
was done in \cite{BO} for $V=M(1)$, in the case
when $u_i=v_i=h(-1){\bf 1}$, $i=1,...,n$, i.e.,
\begin{equation} \label{gen7}
{\cal F}(t_1,...,t_n,q)=
\eta(q) {\rm tr}|_{M(1)} o(Y[h,t_1]h) \cdots o(Y[h,t_n]h)q^{\bar{L}(0)}.
\end{equation}
The approach
in \cite{BO} is different since
it deals with generalized characters and therefore
explicit formulas are easier to obtain. In our approach
we work with the formal delta function which carries
an extra information not present in the generalized character
formulation of (\ref{gen7}) .

\renewcommand{\theequation}{\thesection.\arabic{equation}}
\setcounter{equation}{0}

\section{Relation with vector bundles over elliptic curves}

Let $\theta_{11}(t)$ be as before and $E=\mathbb{C}^\times/q^{\mathbb{Z}}$
be an elliptic curve. Suppose that
${\cal L}(-q^{1/2} t) \in H^0(E, \mathcal{O}_E^*)$ is a
holomorphic line bundle over $E$.
Then ${\cal L}(-q^{1/2} t)$ has no global sections (the degree is -1).
However $F(t)=\frac{1}{\theta_{11}(t)}$ is a multi--valued section of
${\cal L}(-q^{1/2}t)$.

Suppose $n=2$. Let $t_1$ be the local coordinate on $E$ and $t_2$ fixed.
If we set
$$F_1(t_1)=G(t_1,t_2)$$
and
$$F_2(t_1)=G(t_1 t_2),$$
where $G(t_1,t_2)$ is Bloch--Okounkov $2$--point function
associated to a pair of fermions, then \bea &&
F_1(qt_1)=-q^{1/2}t_1 t_2 F_1(t_1)+q^{1/2}t_1t_2 F_2(t_2) \nn &&
F_2(qt_1)=-q^{1/2}t_1 t_2 F_2(t_1). \eea Therefore a pair
$(F_1(t_1),F_2 (t_2))$ can be though as a multi--valued 
section of the vector bundle ${\cal V}$ over $E$ (of
the rank $2$) associated with the conjugacy class (which is enough
to specify according to the classification theory) of the matrix

$$\left ( \begin{array}{cc} -q^{1/2}t_1t_2 & q^{1/2}t_1 t_2 \\ 0 &
-q^{1/2}t_1 t_2 \end{array}
\right ).$$
The latter corresponds to a particular (non--split) extension
of ${\cal L}(-q^{1/2}t_1t_2)$ by itself.

This can be generalized for an arbitrary rank. Fix $t_2,...,t_n$.
For $1 \leq k \leq n$ let \bea \label{very} &&
F_k(t_1)=\sum_{1<i_1<...<i_k \leq n} G(t_1 t_{i_1} \ldots t_{i_k},
\cdots ,\hat{t}_{i_1}, \cdots, \hat{t}_{i_k}, \cdots). \eea Notice
that the right-hand side of (\ref{very}) is invariant with respect
to the symmetric group $S_{n-1}$ (permuting $t_2$,...,$t_n$). Then
$(F_1(t_1),...,F_n(t_1))$ is a multi--valued section
of a rank $n$ vector bundle associated to \bea
&& q^{1/2} \prod_{i=1}^n t_i \left(\begin{array}{cccccc} -1 & 1 & -1  & 1 & ... & ... \\
0 & -1 & { 2 \choose 1 } & -{3 \choose 2} &... & ... \\
0 & 0 & -1 & {3 \choose 1} & ... & ... \\
0 & 0 & 0 & -1 & ... & ... \\
... & ... & ... & ... &.. . & ... \\
... & ... & ... & ... & ... & ...
\end{array} \right).
\eea
\renewcommand{\theequation}{\thesection.\arabic{equation}}
\setcounter{equation}{0}

\section{Conclusion and future work}

\begin{itemize}

\item[(a)]
Notice that one can easily generalize
all results in our paper by considering
more general iterated correlation functions of the form
$$\langle u_{n+1}',X(Y[u_1,y_1-w_1]v_1,e^{w_1}x_1) \cdots  X(Y[u_n,y_n-w_n]v_n,e^{w_n}x_n)
u_{n+1} \rangle,$$
$$\langle u_{n+1}',X(u_1,t_1x_1)X(v_1,s_1 x_1) \cdots X(u_n,t_nx_n)X(v_n,s_n x_n)
u_{n+1} \rangle,$$
and the corresponding $q$--traces.

\item[(b)]
Our considerations in this paper were  mostly at the level
of vertex operator (super)algebras and its modules.
If we closely examine
the $n$--point function $F(t_1,...,t_n)$ associated to
a pair of free charged fermions, we see that we are
actually dealing with certain intertwining operators.
More precisely, the space ${\cal F}_0 \cong M(1)$ is
a vertex operator subalgebra inside the vertex operator superalgebra
${\cal F}$. Therefore operators $X(\psi,x)$ and $X(\psi^*,x)$ do not
act on $M(1)$ (even though the corresponding $2n$--point trace
can be computed for $M(1)$ by using the boson-fermion correspondence).
Thus we are dealing with an non--trivial (even though simple)
intertwining operators of the type
\bea \label{contragrad}
&& {\mathcal{F}_0 \choose \mathcal{F}_{1} \ \mathcal{F}_{-1}} \ \mbox{and} \
{\mathcal{F}_{-1} \choose \mathcal{F}_{-1} \ \mathcal{F}_0},
\eea
where $\mathcal{F}_{\pm 1}$ are $M$--modules.
The fusion algebra for ${\cal F}$, viewed as an intertwining
operator algebra, is isomorphic to
$\mathbb{C}[\mathbb{Z}]$.
Therefore in the future we shall study iterated $2n$--point correlation
functions and corresponding $q$--traces for the abelian intertwining
operator algebras.

\item[(c)]
It is possible to study, over the same lines,
the $n$--point functions associated to certain
$N=2$ vertex operator superalgebras.

\item[(d)] In \cite{M}, motivated by \cite{Bl}, we consider a
large algebra of certain pseudodifferential operators ${\cal
D}^{\pm}_{\infty}$. In the Part III \cite{M3} we study
$n$--point functions twisted with various Dirichlet characters.


\end{itemize}
\renewcommand{\theequation}{\thesection.\arabic{equation}}
\setcounter{equation}{0}

\renewcommand{\theequation}{\thesection.\arabic{equation}}
\setcounter{equation}{0}

\section{Appendix}

\subsection{Appendix A}

In this appendix we prove some elementary results necessary for
solving formal $q$--difference equations for certain $1$-- and
$2$--point functions.

We treat the
following two types
of $q$--difference equations:
\begin{equation} \label{app1}
F(qt)-q^a F(t)=B_1(t)
\end{equation}
and
\begin{equation} \label{app2}
F(qt_1,t_1)+q^{1/2}t_1 t_2 F(t_1,t_2)=B_2(t_1,t_2),
\end{equation}
where $B(t) \in \mathbb{C}[[q^{\pm 1},t^{\pm 1}]]$, $B_2(t_1,t_2) \in
\mathbb{C}[[q^{\pm 1},t_1^{\pm 1},t_2^{\pm 1} ]]$.
We discuss solutions both in $\mathbb{C}[[q^{\pm 1},t_1^{\pm 1},t_2^{\pm 1}]]$
and $\mathbb{C}[[q,t_1^{\pm 1},t_2^{\pm 1}]]$.
This is important since
$$F(qt)-F(t)=D^k_{t} \delta(qt)$$
has two distinguished solutions in
$\mathbb{C}[[q^{\pm 1},t^{\pm 1}]].$
In general (\ref{app1}) may not have any solution.
For example
$$F(qt_1)-F(t)=\delta \left(\frac{t}{q}\right),$$
does not have solution in $\mathbb{C}[[q,t_1^{\pm 1},t_2^{\pm 1}]].$
\begin{remark}
{\em In most of the literature one does not consider formal solutions
of $q$--difference
equations but rather solutions in some analytic spaces (like
holomorphic functions in the disk, punctured disk, etc.).
We work often with series with no convergence (like delta
functions) thus our approach has formal--distribution theoretical flavor.}
\end{remark}

\begin{lemma} \label{appl1}
Let $n \in \mathbb{Z}$.
Then all solutions of
$$F(qt)-q^nF(t)=0,$$
in $\mathbb{C}[[q,t^{\pm 1}]]$
are of the form $t^nf(q)$ for some $f(q) \in \mathbb{C}[[q]]$.
\end{lemma}

\begin{lemma} \label{appl2}
All solutions of
$$F(qt)+q^{1/2}tF(t)=0,$$
in $t^{1/2}\mathbb{C}[[q^{ \pm 1},t^{\pm 1}]]$ are
of the form
\be
f(q) \sum_{n \in \mathbb{Z}} (-1)^n q^{-\frac{n(n+1)}{2}}t^{n+
\frac{1}{2}},
\ee
for some $f(q) \in \mathbb{C}[[q^{\pm 1}]]$. In particular
the only solution inside  $t^{1/2} \mathbb{C} [[q,t^{\pm 1}]]$ is the
trivial one.
\end{lemma}
{\em Proof:}
Since $F(t)=\sum_{n \in \mathbb{Z}} f_{n+1/2}(q) t^{n + \frac{1}{2}},$
it follows that
$$q^n f_{n+1/2}(q)=-f_{n-1/2}(q).$$
Hence, $f_{1/2}(q)=-f_{-1/2}(q)$ and $f_{1/2}(q)$ uniquely determine
$f_{n+1/2}(q)$.

\begin{lemma} \label{appl3}
All solutions of
$$F(qt_1,t_2)-F(t_1,t_2)=0,$$
inside $\mathbb{C}[[q^{\pm 1},t_1^{\pm 1},t_2^{\pm 1}]]$
that satisfy $F(t_1,t_2)=F(t_2,t_1)$ are contained in
$\mathbb{C}[[q^{\pm 1}]]$.
\end{lemma}

\begin{lemma} \label{appl4}
All solutions of
\bea \label{app4}
&& F(qt_1,t_2)+q^{1/2} t_1 t_2 F(t_1,t_2)=0,
\eea
in $t_1^{1/2} t_2^{1/2} \mathbb{C}[[q^{\pm 1},t_1^{\pm 1},t_2^{\pm 1}]]$
that satisfy
\be \label{appsymm}
F(t_1,t_2)=F(t_2,t_1)
\ee
are of the form
$$f(q) \sum_{n \in \mathbb{Z}} (-1)^n q^{-\frac{n(n+1)}{2}}(t_1 t_2)^{n+
\frac{1}{2}}.$$
In particular there is no nontrivial solution inside
$\mathbb{C}[[q,t_1^{\pm 1},t_2^{\pm 1}]]$.
\end{lemma}
{\em Proof:}
Let $F(t_1,t_2)=\sum_{n \in \mathbb{Z}+\frac{1}{2}} t_1^n f_n(t_2,q)$. By
comparing coefficients of $t_i^n$ in (\ref{app4}) we get
$$q^n f_n(t_2,q)+ q^{1/2} t_2 f_{n-1}(t_2,q)=0,$$
On the other hand (\ref{appsymm}) gives us
$$f_n(qt_2,q)+q^{1/2} t_2 f_{n-1}(t_2,q).$$
Combined
\begin{equation} \label{app5}
f_n(qt_2,q)-q^n f_n(t_2,q)=0,
\end{equation}
for every $n \in \mathbb{Z}+\frac{1}{2} $.
Because of Lemma \ref{appl1},
all solutions of (\ref{app5}), inside $t_2^{1/2} \mathbb{C}[[q,t_2^{\pm 1}]]$,
are of the form $t_2^{n} f(q)$, for some $f(q) \in \mathbb{C}[[q^{\pm
1}]]$.
Hence $F(t_1,t_2)=\sum_{n \in \mathbb{Z}+\frac{1}{2}} (t_1t_2) ^n
f_n(q)$.
Now apply Lemma \ref{appl2}.
\epfv

\subsection{Appendix B}

Here we give a different proof of the main technical result used by Zhu
---the recursion formula (cf. Proposition 4.3.4 of \cite{Zh2}) and its consequences.

Instead of working with $P_{k+1}$--functions and their analytic properties (\ref{weier})
we prefer the use formal variable ``all the way'' and then, at the end, recognize
certain formal series as Laurent expansions (around $y=0$) of Weierstrass' functions.
In other words, unlike Zhu's proof, our proof is completely formal---and then at the end---one can ``turn on''
the complex variables

First we (re)prove the following result from \cite{Zh1} (Proposition 4.3.2):
\begin{proposition} \label{431}
\bea \label{432}
&& {\rm tr}|_M X(u_1,x_1)X(u_2,x_2)q^{L(0)}={\rm tr}|_M o(u_1)o(u_2) q^{L(0)}
\nn
&& + \sum_{m \geq 0} {P}_{m+1}(\frac{x_2}{x_1},q)
X(u_1[m]u_2,x_2)q^{L(0)}.
\eea
\end{proposition}
{\em Proof:}
From \cite{Le2},
$$[X(u_1,x_1),X(u_2,x_2)]={\rm Res}_{y} \delta
\left(\frac{e^yx_2}{x_1}\right) X(Y[u_1,y]u_2,x_2).
$$
Now from the property of the trace ${\rm tr}_M (ABC)={\rm tr}|_M (BCA)$
we obtain
\bea \label{z1}
&& (1-q^{-D_{x_1}})X(u_1,x_1)X(u_2,x_2)q^{L(0)}= \nn
&& \sum_{m \geq 0} \frac{D^m}{m!} \delta
\left(\frac{x_2}{x_1}\right) X(u_1[m]u_2,x_2).
\eea
Now, by formally ``inverting'' the operator $1-q^{-D_{x_1}}$ inside
$\mathbb{C}[[q]]$ we obtain the formula (\ref{432}).
\epfv

Let us recall Weierstrass function $\wp_k(z,\tau)$ defined before.
These functions have a Laurent expansion near $z=0$ ($k=1,2,...$) :
\begin{equation} \label{433}
{\wp}_k(z,\tau)=\frac{1}{z^k}+(-1)^k \sum_{n \geq 1} {2n+1 \choose k-1} G_{2n+2}(\tau)z^{2n+2-k}.
\end{equation}
Notice that (\ref{433}) can be considered formally as an element of
$$\mathbb{C}[[z,z^{-1},q]].$$

Let us define related function:
\bea
\bar{\wp}_k(z,q)= \left \{ \begin{array}{cc}
& \tilde{\wp}_2(z,q)+\tilde{G}_2(q) \ {\rm for} \ {\rm for} \ k=2  \\
& \tilde{\wp}_k(z,q), \ k \geq 3  \end{array} \right.
\eea

\begin{theorem} \label{430}
\bea
&& {\rm tr}|_M X(Y[u,y]v,x)=\nn
&& \sum_{m \geq 1} \bar{\wp}_{m+1}(y) {\rm tr}|_M X(u[m]v,x) q^{L(0)}+
{\rm tr}|_M o(u)o(v)q^{L(0)}.
\eea
In particular, if $M$ is a $V$--module and $V$ satisfies $C_2$--condition (cf. \cite{Zh1})
then (\ref{430}) is convergent and it 
has a meromorphic continuation to the whole $y$
plane.
\end{theorem}
Here are some consequences (cf. Propositions 4.3.4 and 4.3.5 in \cite{Zh2}):
\begin{corollary} \label{434}
\begin{equation}
{\rm tr}|_M o(u[-1]v)q^{L(0)}={\rm tr}|_M o(u)o(v)q^{L(0)}+
\sum_{m \geq 1} G_{2k}(q) {\rm tr}|_M o(u[2k-1]v) q^{L(0)}.
\end{equation}
\end{corollary}
Combining Corollary \ref{434} and Proposition \ref{431} we
obtain
\begin{corollary}
\bea
&& {\rm tr}|_M X(u,x_1)X(v,x_2)q^{L(0)}= \\
&& = {\rm tr}|_M o(u[-1]v)q^{L(0)}+\frac{1}{2}{\rm tr}|_M
o(u[0]v)q^{L(0)}+\nn
&& {\rm tr}|_M \sum_{m \geq 0} P_{m+1}\left(\frac{x_2}{x_1},q
\right)o(u[m]v)q^{L(0)}-
\sum_{k \geq 1} {\rm tr}|_M G_{2k}(q) o(u[2k-1]v)q^{L(0)}. \nonumber
\eea
\end{corollary}
\begin{corollary} \label{600}
Let $q=e^{2 \pi i \tau}$. The expression (\ref{430}) is a doubly periodic (elliptic) function with respect to
transformations
$$y \mapsto y+2 \pi i, $$
$$y \mapsto y+2 \pi i \tau.$$
\end{corollary}
{\em Proof of Theorem \ref{430} :}
\bea
&& {\rm tr}|_M X(Y[u,y]v,x_2)q^{L(0)}= \nn
&& = \sum_{i \in \mathbb{Z}} {\rm tr}|_M X(u(i)v,x_2)
e^{y {\rm deg}(u)}(e^y-1)^{-i-1} q^{L(0)} \nn
&& = {\rm Res}_{x_1} \sum_{i \in \mathbb{Z}}
x_2^{{\rm deg}(u)+{\rm deg}(v)-i-1} e^{y {\rm deg}(u)}(e^y-1)^{-i-1} \nn
&& {\rm tr}|_M \left((x_1-x_2)^i Y(u,x_1)Y(v,x_2)-(-x_2+x_1)^i Y(v,x_2)Y(u,x_1) \right) q^{L(0)}= \nn
&& = {\rm Res}_{x_1} \sum_{i \in \mathbb{Z}}
e^{y{\rm deg}(u)}(e^y-1)^{-i-1} x_2^{-1}
\left(\frac{x_2}{x_1} \right)^{{\rm deg}(u)} \cdot \nn
&& {\rm tr}|_M \left( x_2^{-i}(x_1-x_2)^i X(u,x_1)X(v,x_2)-
x_2^{-i}(-x_2+x_1)^i X(v,x_2)X(u,x_1) \right) q^{L(0)}\nn
&&= {\rm Res}_{x_1} \sum_{i \in \mathbb{Z}}
e^{y{\rm deg}(u)} (e^y-1)^{-i-1} x_2^{-1}
\left(\frac{x_2}{x_1} \right)^{{\rm deg}(u)} \cdot \nn
&& \biggl( x_2^{-i}(x_1-x_2)^i \sum_{m \geq 0} {\rm tr}|_M
P_{m+1}\left(\frac{x_2}{x_1},q \right) o(u[m]v)q^{{L(0)}} \nn
&& - x_2^{-i}(-x_2+x_1)^i \sum_{m \geq 0} {\rm tr}|_M P_{m+1}
\left(\frac{q x_2}{x_1},q \right)o(u[m]v) q^{L(0)} \biggr) \nn
&& + {\rm Res}_{x_1} \sum_{i \in \mathbb{Z}}
e^{y{\rm deg}(u)} (e^y-1)^{-i-1} x_2^{-1}
\left(\frac{x_2}{x_1} \right)^{{\rm deg}(u)} \cdot \nn
&& \left\{ \biggl( x_2^{-i}(x_1-x_2)^i - x_2^{-i}(-x_2+x_1)^i \biggr)
{\rm tr}|_M o(a)o(b) q^{L(0)} \right \}. \nn
\eea
We introduce a substitution $t=\frac{x_1}{x_2}$.
Then
\bea \label{499}
&& {\rm Res}_{x_1} \sum_{i \in \mathbb{Z}}
e^{y{\rm deg}(u)} (e^y-1)^{-i-1} x_2^{-1}
\left(\frac{x_2}{x_1} \right)^{{\rm deg}(u)} \cdot \nn
&& \left\{ \biggl( x_2^{-i}(x_1-x_2)^i - x_2^{-i}(-x_2+x_1)^i \biggr)
{\rm tr}|_M o(a)o(b) q^{L(0)} \right\}= \nn
&&= {\rm Res}_{t} \sum_{i \in \mathbb{Z}}
e^{y{\rm deg}(u)} (e^y-1)^{-i-1} t^{-{\rm deg}(u)}
\biggl( (t-1)^i -(-1+t)^i \biggr)
{\rm tr}|_M o(a)o(b) q^{L(0)}= \nn
&& = {\rm Res}_{t} \sum_{i \leq -1}
e^{y{\rm deg}(u)} (e^y-1)^{-i-1} t^{-{\rm deg}(u)}(t-1)^i {\rm tr}|_M o(a)o(b)
q^{L(0)} \nn
&& - {\rm Res}_{t} \sum_{i \leq -1} (-1+t)^i
e^{y{\rm deg}(u)} (e^y-1)^{-i-1} t^{-{\rm deg}(u)}
{\rm tr}|_M o(a)o(b) q^{L(0)} \nn
&& = {\rm Res}_{t} \sum_{i \geq 0} e^{y{\rm deg}(u)}
\frac{(e^y-1)^{i}}{(t-1)^{i+1}} t^{-{\rm deg}(u)} {\rm tr}|_M o(a)o(b) q^{L(0)} \nn
&&- {\rm Res}_{t} \sum_{i \geq 0} e^{y{\rm deg}(u)}
\frac{(e^y-1)^{i}}{(-1+t)^{i+1}} t^{-{\rm deg}(u)}  {\rm tr}|_M o(a)o(b) q^{L(0)} \nn
&& = {\rm Res}_{t} \left(  e^{y{\rm deg}(u)} t^{-{\rm
deg}(u)-1 } \frac{1}{1-e^y/t} {\rm tr}|_M o(a)o(b) q^{L(0)} \right)+ \nn
&& + {\rm Res}_{t}  \left( t^{-{\rm deg}(u)} e^{y{\rm deg}(u)-1}
\frac{1}{1-t/e^y} {\rm tr}|_M o(a)o(b) q^{L(0)} \right) \nn
&&= {\rm tr}|_{M} o(a)o(b)q^{L(0)}.
\eea
Also for every $m \geq 1$
\bea \label{500}
&& {\rm Res}_{x_1} \sum_{i \in \mathbb{Z}}
e^{y{\rm deg}(u)} (e^y-1)^{-i-1} x_2^{-1}
\left(\frac{x_2}{x_1} \right)^{{\rm deg}(u)} \cdot \nn
&& \left\{ x_2^{-i}(x_1-x_2)^i 
P_{m+1}\left(\frac{x_2}{x_1},q \right)- x_2^{-i}(-x_2+x_1)^i P_{m+1}
\left(\frac{q x_2}{x_1},q \right) \right\} \nn
&&= {\rm Res}_{t} \biggl\{ \sum_{i \in \mathbb{Z}}
e^{y{\rm deg}(u)} (e^y-1)^{-i-1} x_2^{-1}
t^{-{\rm deg}(u)} \cdot \nn
&& \left\{ (t-1)^i P_{m+1}(t^{-1},q)-(-1+t)^i P_{m+1}(q t^{-1},q) \right\} \biggr\} \nn
&& = {\rm Res}_t  \biggl\{ 
e^{y{\rm deg}(u)} t^{-{\rm deg}(u)-1 }
\frac{1}{1-e^y/t} P_{m}(t^{-1},q) \nn
&&+ t^{-{\rm deg}(u)} e^{y{\rm deg}(u)-1}
\frac{1}{1-t/e^y} P_{m}(q t^{-1},q) \biggr\} + \nn
&&+ {\rm Res}_t \biggl\{ \frac{e^{y {\rm deg}(u)}}{e^y-1} t^{{\rm deg}(u)} \sum_{i \geq 0}
\frac{(t^{-1}-1)^i}{(e^y-1)^i} \frac{D^{m-1}}{(m-1)!} \delta(t) \biggr\}.
\eea
where we used  the fact that
$$P_m(t,q)-P_m(qt,q)=\frac{D^{m-1}}{(m-1)!} \delta(t).$$
If we combine all $P_m$'s into a single generating function we obtain
$$\sum_{m \geq 1} (P_m(t,q)-P_m(qt,q))x^{m-1}= \delta(e^x t).$$
Now
\bea \label{501}
&& {\rm Res}_t \frac{e^{y {\rm deg}(u)}}{e^y-1} t^{{\rm deg}(u)}
\sum_{i=0}^{\infty} \frac{(t^{-1}-1)^i}{(e^y-1)^i} \delta(e^x t)=\nn
&& = {\rm Res}_t \frac{e^{y {\rm deg}(u)}}{e^y-1} e^{-x {\rm deg}(u)} \sum_{i=0}^{\infty}
\frac{(e^x-1)^i}{(e^y-1)^i} \delta(e^x t) \nn
&& = \frac{ e^{(y-x){\rm deg}(u)} }{e^{y-x}-1}.
\eea
In order to evaluate the last sum in (\ref{500}) we have
to extract the coefficient of $x^{m-1}$ inside (\ref{501}).
Clearly
\bea \label{502}
&& {\rm Coeff}_{x^{m-1}} \frac{e^{(y-x){\rm deg}(u)}}{e^{y-x}-1}= \nn
&&= {\rm Coeff}_{x^{m-1}} e^{-x \frac{d}{dy}} \frac{e^{y {\rm deg}(u)}}{e^y-1} \nn
&&= \frac{(-1)^{m-1}}{(m-1)!} \left(\frac{\partial}{\partial y}\right)^{m-1} \sum_{k=0}^{\infty}
\frac{B_{k}({\rm deg}(u))y^{k-1}}{k !}.
\eea
Now, by using the formula
$$P_{m}(t,q)=\frac{1}{(m-1)!} \left( \sum_{n \geq 1}
\frac{n^{m-1}t^n}{1-q^n}+(-1)^m \sum_{n \geq 1} \frac{n^{m-1}q^n
t^{-n}}{1-q^n} \right),$$ 
we obtain 
\bea \label{503} 
&& {\rm Res}_t  \biggl\{e^{y{\rm deg}(u)} t^{-{\rm
deg}(u)-1 } \frac{1}{1-e^y/t}P_{m}(t^{-1},q) + \nn &&+ 
t^{-{\rm deg}(u)} e^{y{\rm deg}(u)-1} \frac{1}{1-t/e^y}
P_{m}(q t^{-1},q) \biggr\} = \nn 
&& = \frac{(-1)^{m-1}}{(m-1)!}
\sum_{n \geq 0} e^{y({\rm deg}(u)+n)} \frac{q^{{\rm
deg}(u)+n}({\rm deg}(u)+n)^{m-1}}{1-q^{{\rm deg}(u)+n}} \nn && +
\frac{1}{(m-1)!} \sum_{n \geq 1} e^{-yn} \frac{n^{m-1}
q^n}{1-q^n}+ \frac{(-1)^m}{(m-1)!} \sum_{n=1}^{{\rm deg}(u)-1}
e^{yn} \frac{n^{m-1}q^n}{1-q^n}. \eea 
From (\ref{502}) and (\ref{503}) it follows ($m \geq 2$): 

\bea \label{504} 
&&{\rm
Res}_{x_1} \sum_{i \in \mathbb{Z}} e^{y{\rm deg}(u)}
(e^y-1)^{-i-1} x_2^{-1} \left(\frac{x_2}{x_1} \right)^{{\rm
deg}(u)} \cdot \nn 
&& \biggl\{ \left( x_2^{-i}(x_1-x_2)^i
P_{m}\left(\frac{x_2}{x_1},q \right) - x_2^{-i}(-x_2+x_1)^i P_{m}
\left(\frac{q x_2}{x_1},q \right) \right) {\rm tr}|_M
o(u[m-1]v)q^{L(0)} \biggr \}= \nn 
&& = \biggl\{
\frac{1}{y^m}+ \frac{(-1)^{m-1}}{(m-1)!}
\sum_{n \geq 0} \frac{B_{m+n}({\rm deg}(u)) y^n}{n!} +
\frac{(-1)^m}{(m-1)!} \sum_{n=1}^{{\rm deg}(u)-1} e^{yn} n^{m-1}
\nn 
&& + \frac{(-1)^m}{(m-1)!} \left(\sum_{n \geq 1} e^{yn}
\frac{n^{m-1}q^n}{1-q^n}+ (-1)^m e^{-yn} \frac{n^{m-1}q^n}{1-q^n} \right)
\biggr \} {\rm tr}|_M o(u[m-1]v)q^{L(0)}= \nn 
&& = \biggl \{
\frac{1}{y^m} + \frac{(-1)^m}{(m-1)!} \sum_{r \geq
0} \left( \frac{B_{m+r}({\rm deg}(u))}{m+r}-\frac{B_{m+r}}{m+r}
\right) \frac{y^r}{r!} + \nn && \frac{(-1)^m}{(m-1)!} \left(
\sum_{r=0}^{\infty} \sum_{n=1}^{\infty} \frac{n^{r+m-1} q^n
y^r}{r!(1-q^n)}+\frac{(-1)^{m+r} n^{r+m-1}q^n y^r}{r!(1-q^n)}
\right) \biggr \} {\rm tr}|_M o(u[m-1]v)q^{L(0)} \nn 
&& = \biggl
\{ \frac{1}{y^m}+ \frac{(-1)^{m}}{(m-1)!} \sum_{r \geq 0, r+m \in
2\mathbb{Z}} \biggl(-\frac{B_{m+r} y^r}{(m+r)r!}+ \nn 
&&
\frac{2}{(r+m-1)!}\sum_{n=1}^{\infty} {r+m-1 \choose m-1}
\frac{n^{m+r-1}q^n}{1-q^n} \biggr) y^r \biggr\} {\rm tr}|_M
o(u[m-1]v)q^{L(0)} \nn 
&&= \biggl\{ \frac{1}{y^m}+
\frac{(-1)^m}{(m-1)!} \sum_{k \geq m/2, k \in \mathbb{N}} {2k-1 \choose m-1}
\left(\frac{-B_{2k}}{(2k)!}+\frac{2}{(2k-1)!} \sum_{n=1}^{\infty}
\frac{n^{2k-1}q^n}{1-q^n} \right) y^{2k-m} \biggr\} \nn && {\rm
tr}|_M o(u[m-1]v)q^{L(0)} \nn && = \biggl( \frac{1}{y^m}+(-1)^m
\sum_{k \geq m/2, k \in \mathbb{N}} {2k-1 \choose m-1}
\tilde{G}_{2k}(q) y^{2k-m} \biggr) {\rm tr}|_M o(u[m-1]v)q^{L(0)}
\nn
&& = \left \{\begin{array}{ccc} & (
 \tilde{{\wp}}_2(y,q)+\tilde{G}_2(q)) {\rm tr}|_M
o(u[1]v)q^{L(0)}, \ {\rm for} \ m=2 \\ & \\
& \tilde{{\wp}}_{m}(y,q) {\rm tr}|_M o(u[m-1]v)q^{L(0)}, \ {\rm
for} \ {\rm for} \ m>2. \end{array} \right. \nn
&&=\bar{{\wp}}_{m}(y,q) {\rm tr}|_M o(u[m-1]v)q^{L(0)}.
\eea
We used a fact that
$$\tilde{G}_{2k}(\tau)=-\frac{B_{2k}}{(2k)!}+\frac{2}{(2k-1)!}
\sum_{n=1}^{\infty} \frac{n^{2k-1}q^n}{1-q^n}.$$
Because of
\begin{equation} \label{505}
{\rm tr}|_M o(u[0]v) q^{L(0)}=0,
\end{equation}
formulas (\ref{499}), (\ref{504}) and (\ref{505}) imply
\bea
&& {\rm tr}|_M X(Y[u,y]v,x)q^{L(0)}=\nn
&& \sum_{m \geq 1}\bar{\wp}_{m+1}(y,q) {\rm tr}|_M X(u[m]v,x) q^{L(0)}+
{\rm tr}|_M o(u)o(v)q^{L(0)}.
\eea
\noindent {\em Proof of Corollary \ref{600}:} \\
If we let $q=e^{2 \pi i \tau}$, the ellipticity with respect to
$$y \mapsto y+2\pi i,$$
$$y \mapsto y+2 \pi i \tau$$
follows directly from formulas
(\ref{wenorm}) and (\ref{relation}).
\epfv
}
\bibliography{part2}
\bibliographystyle{plain}

\end{document}